\documentclass[11pt]{amsart}
\usepackage{hyperref,amsmath,amsthm,amsfonts,amscd,flafter,epsf,amssymb}
\usepackage{epsfig}
\usepackage{amsfonts}
\usepackage{amssymb}
\usepackage{amsmath}
\usepackage{amsthm}
\usepackage{latexsym}
\usepackage{amscd}
\usepackage{flafter}
\usepackage{epsf}

\newcommand{\ai}{\alpha}
\newcommand{\bi}{\beta}

\newcommand{\lra}{\longrightarrow}
\newcommand{\ra}{\rightarrow}
\newcommand{\Sig}{\Sigma}

\newcommand{\xra}{\xrightarrow}

\newcommand\CF{\widehat{CF}}
\newcommand\CFK{\widehat{CFK}}
\newcommand\HF{\widehat{HF}}

\newcommand\CFL{\widehat{CFL}}
\newcommand\sG{\mathcal G}
\newcommand\sC{\mathcal C}
\newcommand\sD{\mathcal D}

\newcommand\Dual{\mathcal D}
\newcommand\Duality\Dual

\newcommand\RelSpinC{\underline{\SpinC}}
\newcommand\relspinc{\underline{\spinc}}

\newcommand\x{\mathbf x}
\newcommand\w{\mathbf w}
\newcommand\z{\mathbf z}

\newcommand\y{\mathbf y}

\newcommand\Hbar{\overline H}

\newcommand\zZ{\mathcal Z}

\newcommand\sA{\mathcal A}
\newcommand\sB{\mathcal B}

\newcommand\sF{\mathcal F}
\newcommand\sM{\mathcal M}

\newcommand\sZ{\mathcal Z}

\newcommand\ModSphere{\ModFlow\left({\mathbb S}\longrightarrow
\Sym^{g-1}(\Sigma_{1})\times \Sym^2(\Sigma_{2})\right)}
\newcommand\ModSpheres\ModSphere

\newcommand\UnparModSp{\widehat \ModSp}
\newcommand\UnparModFlow\UnparModSp

\newcommand\PD{\mathrm{PD}}

\newcommand{\spinc}{\mathfrak s}

\newcommand{\spinct}{\mathfrak t}

\newcommand\ModMaps{\mathcal M}
\newcommand\ModSp\ModMaps

\newcommand\Ta{{\mathbb T}_{\alpha}}
\newcommand\Tb{{\mathbb T}_{\beta}}
\newcommand\Tc{{\mathbb T}_{\gamma}}

\newcommand\alphas{\mbox{\boldmath$\alpha$}}
\newcommand\thetas{\mbox{\boldmath$\theta$}}

\newcommand\betas{\mbox{\boldmath$\beta$}}
\newcommand\gammas{\mbox{\boldmath$\gamma$}}
\newcommand\deltas{\mbox{\boldmath$\delta$}}

\newcommand\spincrel\relspinc

\newtheorem{thm}{Theorem}[section]
\newtheorem{prop}[thm]{Proposition}

\def\endproof{\relax\ifmmode\expandafter\endproofmath\else
  \unskip\nobreak\hfil\penalty50\hskip.75em\hbox{}\nobreak\hfil\bull
  {\parfillskip=0pt \finalhyphendemerits=0 \bigbreak}\fi}
\def\endproofmath$${\eqno\bull$$\bigbreak}
\def\bull{\vbox{\hrule\hbox{\vrule\kern3pt\vbox{\kern6pt}\kern3pt\vrule}\hrule}}

\newcommand{\R}{\mathbb{R}}

\newcommand{\Z}{\mathbb{Z}}

\newcommand{\ModSWfour}{\mathcal{M}}
\newcommand{\ModFlow}{\ModSWfour}

\newcommand{\SpinC}{{\mathrm{Spin}}^c}

\newcommand\abuts\Rightarrow
\newcommand\Sym{\mathrm{Sym}}

\title{Filtration of Heegaard Floer homology and gluing formulas}
\author{Eaman Eftekhary}

\begin{document}

\begin{abstract}
We introduce an extra filtration on the complex
$\widehat{CFK}(Y,K)$ associated with a null-homologous knot $K$
inside the three-manifold $Y$, denoted by $\CFK_\bullet(Y,K)$,
with $\bullet \in\{-,0,+\}$. This filtration will present the
longitude theory $\widehat{CFL}(Y,K)$ as a subcomplex of
$\widehat{CFK}(Y,K)$. The surgery exact sequences respect this
filtration. Besides some basic properties of these filtered
complexes,
 we derive a formula for $\CFK$ of the knot $(Y,K)$ obtained by gluing the
knot complements $Y_1\setminus \text{nd}(K_1)$ and $Y_2 \setminus \text{nd}(K_2)$.
We will also compute the filtered complex $\CFK_\bullet(S^3,K)$ for an alternating
knot $K$.
\end{abstract}
\maketitle
\section{Introduction}
The introduction of the knot Floer homology by Ozsv\'ath and
Szab\'o (\cite{OS-knot}) and independently by Rasmussen
(\cite{Ras2}) has brought an amazingly powerful technology package
to classical knot theory. These are invariants of the null
homologous knots in three-manifolds, which are constructed
following the general procedure introduced by Ozsv\'ath and
Szab\'o in the context of three-manifolds (\cite{OS-3m1,
OS-3m2}).\\

The knot Floer homology assigned to a knot $(Y,K)$ comes as a
package of filtered chain complexes, the main object being
$CFK^\infty(Y,K)$. There are easier filtered packages which may be
realized as the subcomplexes, quotient complexes, etc., of
$CFK^\infty(Y,K)$ with the help of $\Z \oplus \Z$ filtration
induced by the
knot $K$.\\

In this paper we will be dealing with the filtered chain complex $\CFK(Y,K)$, which is
somewhat easier to compute in general. Suppose that the Heegaard diagram
$$(\Sig,\alphas,\betas_0\cup \{m\},z\}$$
is chosen so that $(\Sig,\alphas,\betas_0)$ represents a Heegaard
diagram for the knot complement $Y\setminus \text{nd}(K)$, which
is completed to a Heegaard diagram for $Y$ by adding the special
$\beta$-curve $m$. Here $\alphas=\{\alpha_1,...,\alpha_g\}$ and
$\betas=\betas_0\cup \{m\}=\{m,\beta_2,...,\beta_g\}$ are two
$g$-tuples of linearly independent simple closed curves on the
genus $g$ surface $\Sig$. Furthermore, assume that $z$ is a marked
point on $m \setminus \alphas$. The points $u$ and $w$ are chosen
on the two sides of the curve $m$ and very close to the marked
point $z$. The diagram should be chosen to be admissible in the
sense of \cite{OS-3m1,OS-knot}, a condition that we will drop from
this exposition. Let $\CFK(Y,K)$ be the complex freely generated
by the intersections of the two tori
$$\Ta=\alpha_1 \times ...\times \alpha_g,\text{ and }
\Tb=m\times \beta_2\times ... \times \beta_g,$$ in the symplectic
manifold $\Sym^g(\Sig)=\frac{\Sig \times ... \times \Sig}{S_g}$,
which is the $g$-th symmetric product of the surface $\Sig$.\\

The differential map of the complex $\CFK(Y,K)$ is defined as follows. For any two
generators $\x,\y \in \Ta \cap \Tb$ let $\pi_2(\x,\y)$ denote the set of homotopy
classes of the disks $\phi:[0,1]\times \R  \ra \Sym^g(\Sig)$ such that $\phi(t,s)$
converges to $\x$ as $s$ goes to $\infty$ and converges to $\y$ as $s$ goes to
$-\infty$, and such that $\phi(0,s)\in \Ta$ and $\phi(1,s)\in \Tb$. For any element
$\phi \in \pi_2(\x,\y)$ let $\sM(\phi)$ denote the moduli space of holomorphic
representatives of $\phi$ with respect to the standard complex structure on the
complex plane, and a generic one parameter family of complex structures
$\{J_t\}_{t\in[0,1]}$ on $\Sym^g(\Sig)$. Let $n_u(\phi)$ and $n_w(\phi)$ denote the
algebraic intersection numbers of the homotopy class $\phi$ with the codimension
$2$ submanifolds $\{u\}\times \Sym^{g-1}(\Sig)$ and $\{w\}\times \Sym^{g-1}(\Sig)$
respectively. Then define

\begin{displaymath}
\partial [\x]=\sum_{\y \in \Ta\cap\Tb}\sum_{
\substack{\phi \in \pi_2(\x,\y)\\ \mu(\phi)=1\\ n_u(\phi)=0}}\#(\frac{\sM(\phi)}
{\sim_\R}).[y],
\end{displaymath}
 where $\sim_\R$ represents the equivalence relation induced by the $\R$-action on
the moduli space $\sM(\phi)$, and $\mu(\phi)$ is the Maslov index associated with
the homotopy class $\phi$.\\

The complex $\CFK(Y,K)$ is filtered using the marked point $w$. In
fact, we may assign a $\SpinC$ structure
$\spinc_u(\x)=\spinc_w(\x)$ to any generator $\x$ of the complex
using either of the marked points $u$ or $w$, which is an element
of $\RelSpinC(Y,K)=\SpinC(Y_0(K))\simeq \Z \oplus \SpinC(Y)$, such
that the boundary map
$\partial$ does not increase the $\Z$-factor of the generators.\\

In this paper, very often
we may refer by $\CFK(Y,K)$ to the above complex with a different
boundary map, which consists only of the part of $\partial$ which does not change
the $\Z$-component of the associated $\SpinC$ structure of the generators in $\CFK(Y,K)$.
Usually, this is the case unless it is clear from the context that we are considering the
whole boundary map.\\

The complex $\CFK(Y,K)$ behaves quite well under the connected sum (\cite{OS-knot}).
However, there are other ways of obtaining a new knot $(Y,K)$ from two given knots
$(Y_1,K_1)$ and $(Y_2,K_2)$, where there is not much known on the behavior of the
Floer homology group $\CFK(Y,K)$ in terms of $\CFK(Y_i,K_i)$, $i=1,2$.\\

We describe here an extra filtration of the complex $\CFK(Y,K)$ by
the elements of the set $\{-,0,+\}$, where the trivial order
$-<0<+$ is assumed on this set. Correspondingly, for any relative
$\SpinC$ structure
$$\relspinc \in \RelSpinC(Y,K)=\SpinC(Y_0(K))$$
we obtain a sequence of subcomplexes
$$\CFK_-(Y,K,\relspinc) \subset \CFK_0(Y,K,\relspinc) \subset \CFK_+(Y,K\relspinc)
=\CFK(Y,\relspinc).$$
We show furthermore that there is a short exact sequence
$$0 \lra \CFK_0(Y,K,\relspinc) \lra \CFK_+(Y,K,\relspinc) \xra{\partial '}
\CFK_-(Y,K,\relspinc-\PD[m]) \lra 0,$$ where $\PD[m]$ denotes the
Poincar\'e dual of the homology class represented by the meridian
$m$ of the knot $K$ in the three-manifold $Y$. The complex
$\CFK_\bullet(Y,K)$ and the map $\partial'$ may be used to
construct a complex
which in fact has the homotopy type of $\CF(Y)$.\\

Our claim in this paper is that this whole package is a topological invariant of the
knot $(Y,K)$ and the $\SpinC$ structure $\relspinc\in \RelSpinC(Y,K)$.

\begin{thm}
The chain homotopy type of the flag of complexes
$$\widehat{CFK}_{-}(Y,K,\relspinc)\xrightarrow{i_{-}}
\widehat{CFK}_0(Y,K,\relspinc)\xrightarrow{i_0}
\widehat{CFK}_{+}(Y,K,\relspinc)$$ is a topological invariant of
the knot $K$ in the three-manifold $Y$, and the relative $\SpinC$
structure $\relspinc \in \RelSpinC(Y,K)=\SpinC(Y_0(K))$. Moreover
the map $\partial'$
  gives a chain map which makes the following  sequence exact:
$$0\rightarrow \widehat{CFK}_0(Y,K,\relspinc)
\xrightarrow{i_0}  \widehat{CFK}_{+}(Y,K,\relspinc) \xrightarrow {\partial'}
 \widehat{CFK}_{-}(Y,K,\relspinc-\text{PD}[m]) \rightarrow 0.$$
The chain homotopy type of this exact sequence is also a topological invariant
of the knot $K$ and the relative $\SpinC$ structure $\relspinc$.
\end{thm}

This extra filtered complex enjoys some of the properties of the non-filtered complex
$\CFK(Y,K)$. For example it  does not depend on the orientation of the knot $K$,
in the following sense:

\begin{prop}
If $(Y,K)$ is a null-homologous knot and if $(Y,-K)$ denotes the
same knot with the opposite orientation, then for any relative
$\SpinC$ structure $\relspinc \in \RelSpinC(Y,K)$, the filtered
complexes  $\widehat{CFK}(Y,K,\relspinc)$ and
$\widehat{CFK}(Y,-K,\relspinc)$ have the same homotopy type.
Moreover if $\relspinc$ extends a torsion $\SpinC$ structure
$\spinc \in \SpinC(Y)$, then the homotopy equivalence shifts the
absolute Maslov grading by $-2k$, where
$$k=\frac{1}{2}\langle c_1(\relspinc),[\widehat{F}]\rangle,$$
for a capping $\widehat{F}$ of a Seifert surface $F$ for the oriented knot
$K$.  In particular
$$H_d(\widehat{CFK}_{\bullet}(Y,K,\relspinc))\simeq
H_{d-2k}(\widehat{CFK}_{\bullet}(Y,-K,\relspinc)),
$$
for any of the filtration levels $\bullet \in \{-,0,+\}$.
\end{prop}

Here we are denoting the homology of a complex $\sC$ by $H_*(\sC)$, where $*$ denotes
the degree of the elements in $\sC$, when there is such a well-defined degree,
 respected by the differential (boundary map) of the complex $\sC$.\\

Although the definition of these complexes uses special types of
Heegaard diagrams for the knots $(Y,K)$, which are in general
considerably more complicated compared to the Heegaard diagrams
used for the computation of non-filtered complex, it is still
possible to do certain explicit computations, as is illustrated in
this paper. As usual, the alternating knots are the first target
of the Floer homology computations. We prove:

\begin{thm}
Suppose that $(S^3,K)$ is an alternating knot in $S^3$, and let $\Delta_K(t)$ and
$\sigma(K)$ denote the Alexander polynomial and the signature of the knot $K$
respectively.
Then the homotopy type of the filtered chain complex $\CFK_\bullet(S^3,K)$ is
completely determined by $\Delta_K(t)$ and $\sigma(K)$.
\end{thm}

In a similar way, we may introduce a filtered version
$\CFL_\bullet(Y,K)$ of the \emph{longitude Floer homology}
$\CFL(Y,K)$. This complex  was introduced in \cite{longitude} as a
variant of the complex $\CFK(Y,K)$ obtained by interchanging the
role of the meridian and the longitude of $K$ in the construction
of the Floer homology groups. There are many similar properties
that the filtered complexes $\CFL_\bullet(Y,K)$ share with
$\CFK_\bullet(Y,K)$. The computations are also very much related.
In particular we will also have:

\begin{thm}
Suppose that $(S^3,K)$ is an alternating knot, and assume that
$\Delta_K(t)$ and $\sigma(K)$ are as in the above theorem.
Then the homotopy type of the filtered chain complex $\CFL_\bullet(S^3,K)$ is also
completely determined by $\Delta_K(t)$ and $\sigma(K)$.
\end{thm}

These filtered complexes are all graded by the $\SpinC$ structures
$\SpinC(Y_0(K))$. In fact, in \cite{longitude} we mentioned that
in the computation of $\CFL(Y,K)$ the two maps
$$\spinc_u,\spinc_w:\Ta \cap \Tb \lra \RelSpinC(Y,K)$$ differ from each other
by a factor $\PD[m]$, where $m$ is a meridian for $K$ as before, and that one may
either choose to filter $\CFL$ using the map $\spinc_u$, or more invariantly
using the \emph{average} $\SpinC$ structure assignment
$$\x \mapsto \spinc(\x)=\frac{\spinc_u(\x)+\spinc_w(\x)}{2}\in \frac{1}{2}\PD[m]+
\RelSpinC(Y,K).$$ For most of the paper \cite{longitude} this
later choice was the convention. Here in the exposition we choose
to filter the complex using the map $\spinc_u$. Then the two
filtered chain complexes $\CFK_\bullet(Y,K)$ and
$\CFL_\bullet(Y,K)$ are very closely related to each other.
Namely,

\begin{prop}
Suppose that $(Y,K)$ is a null homologous knot. Then for any $\SpinC$ structure
$\relspinc \in \RelSpinC(Y,K)$ we have isomorphism of (non-filtered) chain complexes
$$\CFK_-(Y,K,\relspinc-\PD[m])\simeq \frac{\CFK_+(Y,K,\relspinc)}{\CFK_0
(Y,K,\relspinc)}\simeq \CFL(Y,K,\relspinc),\text{ and}$$
$$\CFL_-(Y,K,\relspinc)\simeq \frac{\CFL_+(Y,K,\relspinc)}{\CFL_0(Y,K,\relspinc)}\simeq
\CFK(Y,K,\relspinc),$$ where all the chain complexes are assumed
to be graded using the $\SpinC$ structures in $\RelSpinC(Y,K)$.
\end{prop}

The techniques used in \cite{longitude} may be used to derive the
following genus formula for the knots in $S^3$. Suppose that
$(S^3,K)$ is a knot, and identify the set of $\SpinC$ structures
$\RelSpinC(S^3,K)=\SpinC(S^3_0(K))$ with $\Z$ by evaluating the
Chern class of a $\SpinC$ structure on a capped Seifert surface
for $K$. For a complex $\sC$ which is graded by $\Z$ (in the sense
that the differential of $\sC$ takes elements in a grading level
$s$ to elements in grading level $s$), define the degree
$d_+(\sC)$ to be given by
$$d_+(\sC)=\text{max}\Big\{s\in \Z \ |\ \ \sC(s)\text{ does not have trivial chain
homotopy type} \Big \},$$
and similarly define
$$d_-(\sC)=\text{min}\Big\{s\in \Z \ |\ \ \sC(s)\text{ does not have trivial chain
homotopy type} \Big \}.$$
Here we are denoting the part of complex $\sC$ in grading level $s$ by $\sC(s)$.
Then we have the following theorem.

\begin{thm}
Suppose that $K$ is a knot in $S^3$ of genus $g$. Then for the complex
$\CFK(K)=\CFK(S^3,K)$ we have
$$d_+(\CFK_-(K))+1=-d_-(\CFK_-(K))=g,$$
$$d_+(\CFK_0(K))=-d_-(\CFK_0(K))=g,\text{ and}$$
$$d_+(\CFK_+(K))=-d_-(\CFK_+(K))=g.$$
Furthermore, for the complex $\CFL(K)=\CFL(S^3,K)$ we have
$$d_+(\CFL_-(K))=-d_-(\CFL_-(K))=g,$$
$$d_+(\CFL_0(K))=-d_-(\CFL_0(K))=g,\text{ and}$$
$$d_+(\CFL_+(K))=1-d_-(\CFL_+(K))=g.$$
\end{thm}

As we mentioned before, for many gluing constructions, it is not known how to relate
the Floer homology of the final object to the Floer homology of the building blocks
of the construction. The main type of construction we have in mind is gluing of
three manifold along certain null-homologous knots inside them.\\

Suppose that $(Y_1,K_1)$ and $(Y_2,K_2)$ are two null-homologous
knots and consider the knot complements $W_i=Y_i \setminus
\text{nd}(K_i)$. There are two distinguished curves on the torus
boundary of each $W_i$, $i=1,2$. One of them is the meridian $m_i$
of the knot $K_i$ in $Y_i$, and the other one is a longitude $l_i$
for $K_i$, so that it gives the three manifold $(Y_i)_0(K_i)$,
obtained from $Y_i$ by a zero surgery on $K_i$. The two curves
$(m_i,l_i)$, $i=1,2$, determine a framing of the boundary of the
three-manifold $W_i$. One may glue $W_1$ to $W_2$ along their
torus boundary in many ways. Of special interest to us are the
following two special cases. The first case is when the curve
$m_1$ is glued to $m_2$ and the curve $l_1$ is glued to $l_2$
under the above identification of the boundaries of $W_1$ and
$W_2$. The identification of $l_1$ with $l_2$, in fact, identifies
the two knots $K_1$ and $K_2$ to give a knot $K$ in the resulting
three-manifold $Y$. The result of this construction
 will be denoted by
$(Y,K)=(Y_1,K_1)\|(Y_2,K_2)$. We may also choose the
identification of the boundaries so that $m_1$ is identified with
$l_2$ and $m_2$ is identified with $l_1$. A parallel copy of $K_1$
gives a knot $K$ in the resulting three-manifold $Y$, and we will
denote the result of this construction by $(Y,K)=(Y_1,K_1)\perp
(Y_2,K_2)$. The Heegaard Floer homology groups (non-filtered
versions) $\CFK(Y,K)$ and $\CFL(Y,K)$ may be computed in terms of
filtered Floer homologies of $(Y_i,K_i)$ for these constructions.
In particular, we prove the following two parallel theorems:

\begin{thm}
Suppose that $(Y_1,K_1)$ and $(Y_2,K_2)$ are two null-homologous knots and that
$(Y,K)=(Y_1,K_1)\|(Y_2,K_2)$ as above. Then for any
relative $\SpinC$ structure
$$\relspinc_1\#\relspinc_2 \in \RelSpinC(Y,K)=\RelSpinC(Y_1,K_1)\oplus
\RelSpinC(Y_2,K_2),$$
the
Heegaard Floer homology $\CFK(Y,K)$ in the $\SpinC$ structure $\relspinc_1\#\relspinc_2$
will be given as the quotient complex
\begin{displaymath}
\sC(\relspinc_1\#\relspinc_2)=\frac{\Big [
\frac{\CFL(Y_1,K_1,\relspinc_1)\otimes \CFL(Y_2,K_2,\relspinc_2)}
{\CFL_0(Y_1,K_1,\relspinc_1)\otimes
\CFL_0(Y_2,K_2,\relspinc_2)}\Big ] }{\Large{\sim}},
\end{displaymath}
where the equivalence relation $\sim$ is induced by the isomorphism of the
subcomplexes
\begin{displaymath}
\begin{split}
\rho: \CFL_-(Y_1,K_1,\relspinc_1)&\otimes \frac{\CFL(Y_2,K_2,\relspinc_2)}
{\CFL_0(Y_2,K_2,\relspinc_2)}\\
&\lra
\frac{\CFL(Y_1,K_1,\relspinc_1)}{\CFL_0(Y_1,K_1,\relspinc_1)}\otimes
\CFL_-(Y_2,K_2,\relspinc_2).
\end{split}
\end{displaymath}
Here $\rho$ is induced by the isomorphisms $\CFL_-\simeq
\frac{\CFL}{\CFL_0}$. Moreover, the groups $\CFL(Y,K)$ are given
by
$$\CFL(Y,K,\relspinc_1\#\relspinc_2)\simeq \CFL(Y_1,K_1,\relspinc_1)\otimes
\CFL(Y_2,K_2,\relspinc_2).$$
\end{thm}

\begin{thm}
Suppose that $(Y,K)=(Y_1,K_1)\perp (Y_2,K_2)$ as above. Then
$$\RelSpinC(Y,K)\oplus \Z=\RelSpinC(Y_1,K_1) \oplus \RelSpinC(Y_2,K_2),$$
with $\Z$ being generated by $\text{PD}[m_2]$. Furthermore, for any
relative $\SpinC$ structure $\relspinc\in \RelSpinC(Y,K)$ we have
\begin{displaymath}
\begin{split}
\CFL(Y,K,\relspinc)&=\bigoplus_{\substack{\relspinc_i \in \RelSpinC(Y_i,K_i)\\
(\relspinc_1,\relspinc_2) \text{ induce } \relspinc}} \CFL(Y_1,K_1,\relspinc_1)
\otimes \CFK(Y_2,K_2,\relspinc_2),\\
\CFK(Y,K,\relspinc)&=\frac{\Big [\bigoplus_{\substack{\relspinc_i \in \RelSpinC(Y_i,K_i)\\
(\relspinc_1,\relspinc_2) \text{ induce } \relspinc}} \sC(\relspinc_1,\relspinc_2) \Big ]}
{\Large \sim},\\
\sC(\relspinc_1,\relspinc_2)&=\frac{\CFL(Y_1,K_1,\relspinc_1)\otimes
\CFK(Y_2,K_2,\relspinc_2)}{\CFL_0(Y_1,K_1,\relspinc_1)\otimes
\CFK_0(Y_2,K_2,\relspinc_2)},\\
\CFL_-(Y_1,K_1,\relspinc_1)&\otimes
\frac{\CFK(Y_2,K_2,\relspinc_2)}{\CFK_0(Y_2,K_2,\relspinc_2)}{\Large \sim}\\
&\frac{\CFL(Y_1,K_1,\relspinc_1)}{\CFL_0(Y_1,K_1,\relspinc_1)}\otimes
\CFK_-(Y_2,K_2,\relspinc_2-\text{PD}[m_2]).
\end{split}
\end{displaymath}
As before the equivalence relation $\sim$ is induced by the connecting isomorphisms
$\CFK_-\simeq \frac{\CFK}{\CFK_0}$ and $\CFL_-\simeq \frac{\CFL}{\CFL_0}$.
\end{thm}

These gluing formulas are in fact our main justification, besides
their nice properties, for definition of these extra filtered
complexes. Some of the computation technology introduced in
\cite{OS-knot,OS-alternating} naturally generalizes to these extra
filtered versions of our complexes. One thing we would like to
emphasize here, is the surgery long exact sequence. If $(Y,K)$
denotes a knot in a three-manifold $Y$, and if $\gamma$ is another
knot which has zero linking number with $K$, we may extend the
surgery short exact sequence of Ozsv\'ath and Szab\'o to the
context of filtered chain complexes to obtain the following long
exact sequence formula in the level of homology.

\begin{thm}
Let $(Y,K)$ be a knot and let $\gamma$ be a framed knot in $Y$ which
is disjoint from $K$ and has zero linking number with it. Then for a correct
choice of the Seifert surface, and for each integer
$m\in \mathbb{Z}$ we obtain the exact sequences:
\begin{displaymath}
\begin{split}
...\ra H_*({\CFK}_{\bullet}(Y_{-1}(\gamma),K,m))
&\xra{f^1_{\bullet}}
 H_*({\CFK}_{\bullet}(Y_{0}(\gamma),K,m)) \xra{f^2_{\bullet}}\\
&\ra H_*({\CFK}_{\bullet}(Y,K,m)) \xra{f^3_{\bullet}} ...,
\end{split}
\end{displaymath}
where $\bullet \in \{-,0,+\}$.
The maps $f^{\bullet}_1$ and $f^{\bullet}_2$, when the groups are
graded, each will lower the absolute grading by $\frac{1}{2}$,
and $ f^{\bullet}_3$ will not increase the absolute grading.
\end{thm}

{\bf Acknowlegement.} The study of gluing formulas was suggested to 
me by Peter Kronheimer and I would like to thank him for this suggestion.
I would also like to thank Zolt\'an Szab\'o and Jake Rasmussen for interesting 
discussions.

\section{Construction of the filtration}
We begin by considering a special Heegaard diagram for the null-homologous
knot $K$ in the three-manifold $Y$. To obtain such a diagram, start with
a weakly admissible marked Heegaard diagram
$$(\Sig, \alphas,\betas_0\cup\{m\},z)
$$
for $K\subset Y$, as defined by Ozsv\'ath and Szab\'o in
\cite{OS-knot}. As usual $\alphas=\{\ai_1,...,\ai_g\}$ and
$\betas_0=\{\bi_2,...,\bi_g\}$ give the Heegaard diagram for the
knot complement $Y\setminus \text{nd}(K)$ and $m$ represents a
meridian of $K$ in $Y$. Furthermore, the marked point $z$ is
located
on $m$. \\

We may consider a longitude $l$ for the knot $K$ on this surface $\Sigma$
with the property that it is disjoint from the curves in $\betas_0$, cuts
$m$ transversely in a single point and such that the Heegaard diagram
$$(\Sig,\alphas,\betas_0\cup \{l\})
$$
represents the three manifold $Y_0(K)$ obtained from $Y$ by a zero surgery
on $K$. Such a curve $l$ is unique up to handle slides along the curves
in $\betas_0$ (and also isotopies which are disjoint from $\betas$).

We may assume that $l$ cuts the meridian $m$ in the marked point $z$,
without loosing the generality.\\

To get the desired Heegaard diagram, modify this data as follows:
Let $\delta$ be a small oriented arc on $l$, forming a
neighborhood of $z$ in $l$ and disjoint from all the other curves
in $\alphas$ and $\betas_0$. Let $P$ and $Q$ be the start point
and the end point of this arc.  Attach a handle to the surface
$\Sig$ which connects $P$ to $Q$. Using this handle we may replace
the arc $\delta$ on $l$  by a path which avoids $m$, to get a new
curve $\lambda$. The loop $\mu$ which goes once around the handle
cuts $\lambda$ in a single point. Moving $\mu$ by an isotopy, we
may also
create a pair of cancelling intersection points with $m$.\\

Let us denote the new surface by $\Sig '$. If we connect the end
points of $\delta$ to each other using the handle, so that the
resulting closed curve $l'$ only has a single intersection point
with $m$ and a single intersection point with $\mu$, then the
Heegaard diagram

$$(\Sig',\alphas \cup \{\mu\},\betas_0\cup \{\lambda\} \cup \{m\},z)
$$
also represents the same knot $K \subset Y$, and $l'$ represents a longitude
for $K$ in this new diagram (in the above sense).\\

\begin{figure}
\mbox{\vbox{\epsfbox{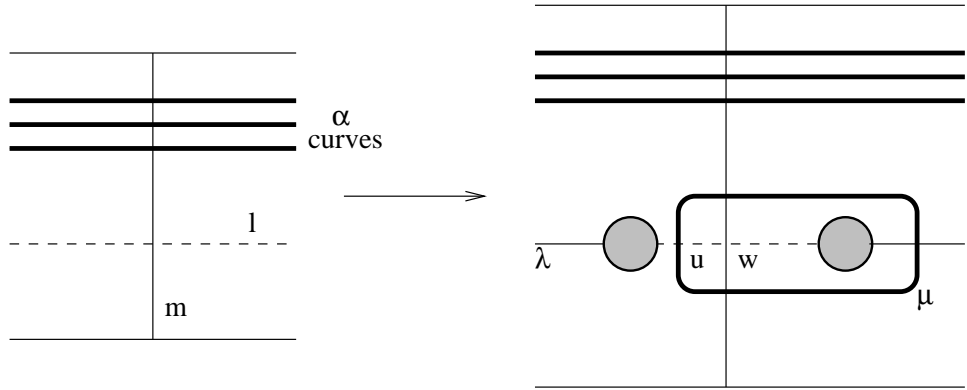}}}
\caption{\label{fig:diagram} {From the given Heegaard diagram on
$\Sig$, where the $\alpha$ curves are represented by the bold
curves, the elements of $\betas_0$ and the meridian $m$ are the
normal lines, and $l$ is the dotted curve, we may obtain a new
Heegaard diagram for the same knot, by adding a handle and the two
curves $\lambda$ and $\mu$. The new longitude $l'$ is the dotted
curve on the right hand side.}}
\end{figure}

Note that given the orientation on $m$ and on $l$, the above construction
may be done in a unique way, say as is suggested by figure~\ref{fig:diagram}.\\

\begin{figure}
\mbox{\vbox{\epsfbox{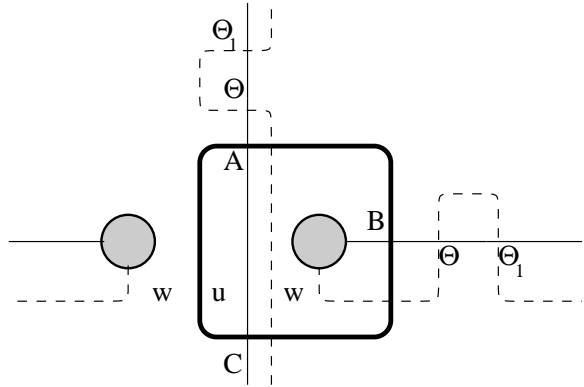}}}
\caption{\label{fig:perturbation}
{When a handle slide happens between the $\beta$ circles in the complement
of the marked point $z$, we may consider a Hamiltonian isotope of the curves
$\lambda$ and $m$ which are denoted by the dotted curve. We assume that
there are a pair of cancelling intersection points between each curve and
its Hamiltonian translate, denoted by $\Theta$ and $\Theta_1$.
}}
\end{figure}

There are three intersection points on the curve $\mu$ which are
named by the letters $A,B$ and $C$ in
figure~\ref{fig:perturbation}. Any generator corresponding to the
Heegaard diagram
$$(\Sig',\alphas \cup \{\mu\},\betas_0\cup \{\lambda\} \cup \{m\},z)
$$
will contain exactly one of these three intersection points, and
accordingly we may partition the generators into three types: type
A, type B and
type C.\\

Put the marked points $u$ and $w$ on the left and on the right hand side of
the point $z\in m$ in order to get a doubly pointed Heegaard diagram, as
in \cite{OS-knot}.\\

Note that any boundary map starting from a generator of type C,
will go to a generator of type C. The reason is that if a homotopy
type $\phi \in \pi_2(\x,\y)$ of disks connecting the generators
$\x$ and $\y$ contains a holomorphic representative, then the
associated domain $\mathcal{D}(\phi)$ will only have positive
coefficients. Since three of the four quadrants around the
intersection point $C$ are forced to have multiplicity zero, the
fourth domain will have a coefficient of $\pm 1$, if $\phi$
connects a generator of type C to a generator of a different type.
Having a holomorphic representative means that this
coefficient is in fact $+1$.\\

But then the differential map associated with this disk will ``end
up'' within the generator of
type C, rather than ``start from'' such a state.\\

A similar argument shows that the boundary of generators of type B
consists of generators of types B and C. As a result we have a
filtration associated with the ordering
$$\text{type C}<\text{type B}<\text{type A}.$$

According to the general properties of filtered chain complexes we
get a filtration of the complex $\widehat{CFK}(Y,K)$ by a flag
$$\widehat{CFK}_{-}(Y,K)\xrightarrow{i_{-}}
\widehat{CFK}_{0}(Y,K)\xrightarrow{i_0}
\widehat{CFK}_{+}(Y,K)=\widehat{CFK}(Y,K).$$

Here $\widehat{CFK}_{-}(Y,K)$ denotes the subcomplex generated by
the generators of type C, $\widehat{CFK}_{0}(Y,K)$ consists of
generators of types B and C, and $\widehat{CFK}_{+}(Y,K)$ consists
of all the generators (of types A,B and C). The maps $i_{-}$ and
$i_0$ are embeddings of the subcomplexes. The boundary maps are
the restrictions of the standard boundary map $\partial_+$ of
$\widehat{CFK}(Y,K)$ to the
corresponding subcomplex.\\

Before we state the invariance of this filtration, let us
introduce one extra
map $\partial'$  on the complex.\\

If we allow the connecting disks to have a nonzero coefficient at
the point $u$, i.e. if we allow $n_u(\phi)\neq 0$, then the only
possible disks between two generators with $n_u(\phi)\neq 0$ which
also have holomorphic representatives will go from  a generator of
type A to a generator of type C. Furthermore, for such a disk,
$n_u(\phi)=1$ and the other regions which have $A$ or $C$ as a
corner will get a zero coefficient
in the domain $\mathcal{D}(\phi)$ of $\phi$.\\

If the two generators agree, except for the intersection points
$A$ and $C$, then there is a unique holomorphic disk connecting them.\\

As a result of this discussion we see that the map $\partial_1$,
which counts such disks, is trivial except on the generators of
type $A$, and that the image of $\partial_1$ lies
in $\widehat{CFK}_{-}(Y,K)$.\\

It is important to note that the complex $\widehat{CFK}(Y,K)$
equipped with $\partial=\partial_++\partial_1$ gives the chain
complex which evaluates the Floer homology of the three-manifold
$Y$ (the hat theory), where $\partial_+$ is the differential of
$\CFK_+(Y,K)$. We may combine the map $\partial_1$ with the sign
map, assigning the values $\pm 1$ to the generators of the complex
according to their absolute
$\frac{\Z}{2\Z}$-grading to get a chain map $\partial'=\epsilon \partial_1$.\\

We also make the remark that as usual the $\SpinC$ structures of
$Y_0(K)$ will filter the complexes:
$$\widehat{CFK}_{-}(Y,K,\relspinc)\xrightarrow{i_{-}}
\widehat{CFK}_0(Y,K,\relspinc)\xrightarrow{i_0}
\widehat{CFK}_{+}(Y,K,\relspinc)=\widehat{CFK}(Y,K,\relspinc).$$
Here $\widehat{CFK}(Y,K)=\bigoplus_{\relspinc \in \SpinC(Y_0(K))}
\widehat{CFK}(Y,K,\relspinc)$, etc..\\

The first result of this paper is the following:

\begin{thm}
The chain  homotopy type of the flag of complexes
$$\widehat{CFK}_{-}(Y,K,\relspinc)\xrightarrow{i_{-}}
\widehat{CFK}_0(Y,K,\relspinc)\xrightarrow{i_0}
\widehat{CFK}_{+}(Y,K,\relspinc)$$ is a topological invariant of the knot
$K$ in the three-manifold $Y$, and the relative $\SpinC$ structure
$\relspinc \in \RelSpinC(Y)=\SpinC(Y_0(K))$. Moreover the map $\partial'$
 satisfies
$\partial' \circ \partial_+-\partial_+ \circ \partial'=0$ which
gives the following  exact sequence:
$$0\rightarrow \widehat{CFK}_0(Y,K,\relspinc)
\xrightarrow{i_0}  \widehat{CFK}_{+}(Y,K,\relspinc) \xrightarrow {\partial'}
 \widehat{CFK}_{-}(Y,K,\relspinc-\text{PD}[m]) \rightarrow 0.$$
The chain homotopy type of this exact sequence is also a
topological invariant of the knot $K$ and the relative $\SpinC$
structure $\relspinc$.
\end{thm}
\begin{proof}
The independence from the almost complex structure is easy and standard.
We have to show that different choices of the Heegaard diagram do not
change the homotopy types of the the above sequences.\\

To do so, we should show that:\\

1) The final choice of $l$, up to a handle slide does not change the
homotopy types.

2) Handle slides and isotopies of $\alphas$, away from the base
point $z$, preserve the homotopy types.

3) Handle slides of $\betas_0$ and $m$ along $\betas_0$ keep the
homotopy type invariant.

4) The same is true for the isotopies of $\betas_0$ and $m$.

5) The homotopy type is not changed in the process of handle addition.\\

For the first claim, assume that $l^1$ is obtained from $l$ by a
handle slide along one of the curves in $\betas_0$, say $\beta_2$.
We implicitly assume that $l^1$ is moved by a hamiltonian isotopy
so that a pair of cancelling transverse intersections are created
between $l$ and $l^1$. Correspondingly introduce $\lambda^1$,
which is constructed from $l^1$
in the same way that $\lambda$ is constructed from $l$. \\

Let $m^1,\gamma_2,\gamma_3,...,\gamma_g$ be Hamiltonian isotopes
of $m,\beta_2,...,\beta_g$, so that there are  pairs of cancelling
intersection points between $\beta_i$ and $\gamma_i$, and between
$m$ and $m^1$. Then the three sets of curves:
$$\alphas'=\alphas\cup \{\mu\}=\{\alpha_1,...,\alpha_g,\mu\},$$
$$\betas'=\betas \cup \{\lambda\}=\{m,\beta_2,...,\bi_g,\lambda\},
\text{ and}$$
$$\gammas'=\gammas \cup \{\lambda^1\}=\{m^1,\gamma_2,...,\gamma_g,\lambda^1\}$$
will form a Heegaard triple. $\alphas'$ and $\betas'$ will give
the initial diagram $H_1$, the pair $(\betas',\gammas')$ will give
a Heegaard diagram $H_2$ for the connected sum of $g$ copies of
$S^1\times S^2$, and finally the pair $(\alphas',\gammas')$ will
give the diagram $H_3$ obtained by a handle slide on $H_1$. The
marked points $u$ and $w$ will be fixed for each of the three
Heegaard diagram obtained from the three pairs of $\alphas,\betas$
and $\gammas$. However, on $H_2$, they are both located
in the same region. \\

 The Floer homology of the Heegaard diagram
$$H_2=(\Sig',\betas',\gammas';u=w)$$
gives $\HF(\#^g(S^1\times S^2))$. There is a top generator of this
Floer homology group in the $\SpinC$ structure with trivial first
Chern class, which will be denoted by $\Theta$. The map
$$\mathcal{G}:\widehat{CFK}(Y,K,\relspinc;H_1)\lra
\widehat{CFK}(Y,K,\relspinc;H_3)$$ is defined by
$\mathcal{G}(\x)=\mathcal{F}(\x \otimes \Theta)$ for any generator
$\x \in \Ta \cap \Tb$. Here
$$\Ta=\mu \times \ai_1 \times ...\times \ai_g,$$
$$\Tb=\lambda \times m \times \beta_2 \times ... \times \bi_g, \text{ and}$$
$$\Tc=\lambda^1 \times m^1 \times \gamma_2 \times ... \times \gamma_g,$$
and the map
$$\mathcal{F}:\widehat{CFK}(H_1) \otimes \widehat{CFK}(H_2)
\lra \widehat{CFK}(H_3)$$ is defined for any $\x \otimes \y$ with
$\x \in \Ta \cap \Tb$ and $\y\in \Tb \cap \Tc$ to be

$$\mathcal{F}(\x \otimes \y)=\sum_{\w \in \Ta \cap \Tc}
\sum_{\substack{\phi \in \pi_2(\x,\y,\w)\\
               \mu(\phi)=0\\
               n_u(\phi)=n_w(\phi)=0 }} \#(\mathcal{M}(\phi)).\w.$$
The similar arguments in \cite{OS-knot, OS-3m1, OS-3m2} may be
followed to show that this map is in fact a chain homotopy
equivalence. However, we need to show furthermore, that this map
in fact respects the filtrations of the
complexes $\widehat{CFK}(Y,K;H_1)$ and $\widehat{CFK}(Y,K;H_3)$.\\

It suffices to show that when $\x$ is a generator of type C, then
only generators of type C can appear in the formal sum
$\mathcal{G}(\x)$, and that if $\x$ is of type B, then every
generator that appears in the formal expression $\mathcal{G}(\x)$
is either
of type B or of type C.\\

Note that if the coefficient of $\w \in \Ta\cup \Tc$ is positive in
 $\mathcal{G}(\x)$, then we should have a holomorphic representative for
a homotopy class of disks $\phi \in \pi_2(\x,\Theta,\w)$ with
$n_u(\phi)=n_w(\phi)=0$. If $\x$ is of type C, as we travel
counterclockwise on $\mu$ starting at $C$, before getting to the
intersection point corresponding to $\w$, the coefficient on the
right hand side, is always strictly lower than the one on the left
hand side. Here we are referring to a picture which is locally
illustrated in figure~\ref{fig:perturbation}. If $\w$ is not of
type C, then at some point, we will have a zero coefficient on the
left, which implies that the coefficient on the right is forced to
be
negative, a contradiction.\\

A similar argument shows that for any $\x \in \Ta\cap \Tb$ of type
B, $\mathcal{G}(\x)$ will only consist of nonzero multiples of
generators in $\Ta \cap \Tc$ of types B and C.\\

This proves that $\mathcal{G}$ respects the filtration and that
the restriction to each of the subcomplexes is a homotopy
equivalence as well.\\

To show that the homotopy type of the exact sequences are
preserved, we have to show that for the vertical maps
$\mathcal{G}$ in the diagram

\begin{equation}
\begin{array}{ccccccccc}
0& \ra & \widehat{CFK}_{0}(H_1) &\xra{i_0} &
\widehat{CFK}_{+}(H_1)
&\xra{\partial'} & \widehat{CFK}_{-}(H_1) &\ra & 0\\
&&{\Big\downarrow\vcenter{%
\rlap{$\mathcal{G}$}}}&& {\Big\downarrow\vcenter{%
\rlap{$\mathcal{G}$}}} &&
{\Big\downarrow\vcenter{%
\rlap{$\mathcal{G}$}}} &&\\
0& \ra & \widehat{CFK}_0(H_3) &\xra{i_0} & \widehat{CFK}_{+}(H_3)
&\xra{\partial'} & \widehat{CFK}_{-}(H_3) &\ra & 0,
\end{array}
\end{equation}
the expressions $\mathcal{G}\circ \partial' -\partial' \circ
\mathcal{G}$ and $\mathcal{G}\circ i_B -i_B \circ \mathcal{G}$ are
null homotopic. This is also a typical application of techniques
in \cite{OS-3m2} and
\cite{OS-knot}.\\

For the other cases, the argument continues exactly in the same
way. We will follow the steps showing that the moves of type 2,3,4
and 5 do not change the homotopy type of the complex
$\widehat{CFK}(Y,K)$ as in \cite{OS-knot}, and at each step, we
will also check that the filtration is also preserved by the
homotopy equivalence suggested by the argument of Ozsv\'ath and
Szab\'o in
\cite{OS-knot}. This completes the proof.\\
\end{proof}

It is interesting to note that the homology of the subcomplex
$\widehat{CFK}_{-}(Y,K)$, gives the longitude Floer homology of the knot
$K$ introduced in \cite{longitude}:

\begin{thm}
The homology of the subcomplex $\widehat{CFK}_{-}(Y,K)$ and the quotient
complex
$$\frac{\widehat{CFK}_{+}(Y,K)}{i_0(\widehat{CFK}_0(Y,K))}$$
give the longitude
Floer homology $\widehat{HFL}(Y,K)$. More precisely:
$$H_*(\widehat{CFK}_{-}(Y,K,\relspinc))\simeq \widehat{HFL}(Y,K,\relspinc+
\frac{1}{2}\text{PD}[m]), \text{ and}$$
$$H_*(\frac{\widehat{CFK}_{+}(Y,K,\relspinc)}
{i_0(\widehat{CFK}_0(Y,K,\relspinc))})\simeq
\widehat{HFL}(Y,K,\relspinc-\frac{1}{2}\text{PD}[m]),$$
where PD$[m]$ denotes the Poincar\'e dual of the meridian $m$ of the knot
$K$ in the three-manifold $Y$.
\end{thm}

\begin{proof}
Note that if we replace the meridian $m$ of the diagram shown in
figure~\ref{fig:diagram} by the dotted curve $l'$, we obtain a weakly
admissible diagram for the zero surgery $Y_0(K)$ on the knot $K$. This
diagram may be used for the longitude Floer homology (the hat theory),
as any periodic domain has both positive and negative coefficients.\\

The only curve which cuts $l'$ is the curve $\mu$. It cuts $l'$ in
a single point $z$ and any intersection  of the tori associated
with the tuples $\alphas \cup \{\mu\}$ and $\{l',\lambda\}\cup
\betas_0$ will consist of a $(g+1)$-tuple of points, where one of
the points is forced to be $z$. These $g$ tuples are in one to one
correspondence with the intersection points of $\Ta$ and $\Tb$ of
type $C$. Namely, any $(g+1)$-tuple $\{z,\bullet\}$
is in correspondence with $\{C,\bullet\}$.\\

Since the disks are forced to have zero coefficients at $u$ and
$w$, we conclude that the domains of the disks between the
generators in the longitude theory, are in fact exactly the same
as the domains of
the disks between the corresponding generators of type C.\\

This gives the first isomorphism claimed in the above theorem. The
second
isomorphism is quite similar.\\

The claim about the $\frac{1}{2}\text{PD}[m]$ shifts in the $\SpinC$ structure,
is because of our special convention in \cite{longitude}, on averaging
the $\SpinC$ structures assigned to a generator $\x$ by the maps $\spinc_u$
and $\spinc_w$.\\
\end{proof}

\begin{figure}
\mbox{\vbox{\epsfbox{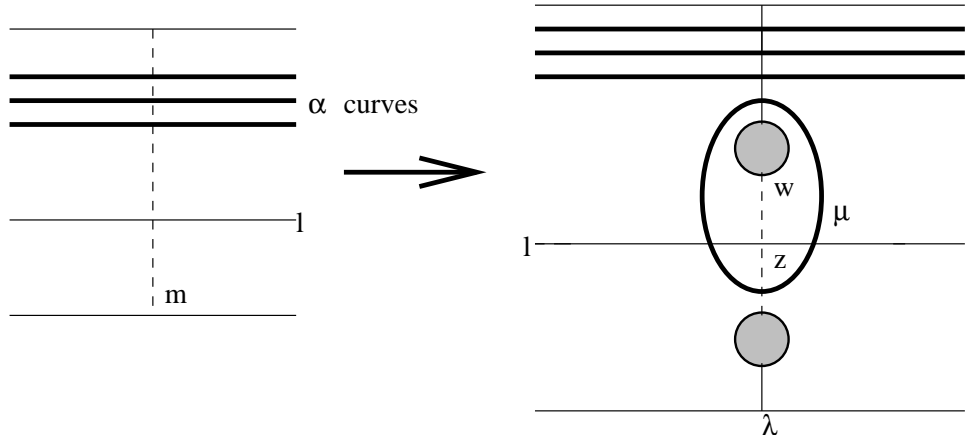}}}
\caption{\label{fig:diagram2} {We may obtain a new Heegaard
diagram for the longitude theory of the knot $K$, representing the
three manifold $Y_0(K)$, by adding a handle and the two curves
$\lambda$ and $\mu$ to the collections $\betas_0\cup \{l\}$ and
$\alphas$ respectively.}}
\end{figure}

There is a second way of getting a filtration, this time on
$\widehat{CFL}(Y,K)$ rather than on $\widehat{CFK}(Y,K)$. Namely,
in the initial diagram of figure~\ref{fig:diagram}, we replace the
role of the meridian $m$ and the longitude $l$ and will continue
with the construction, as is suggested in
figure~\ref{fig:diagram2}. The points $u$ and $w$ will give us two
maps
$$\spinc_u,\spinc_w:\Ta \cap \Tb \lra \RelSpinC(Y,K)=\SpinC(Y_0(K)),$$
which are related by the formula
$$\spinc_u(\x)=\spinc_w(\x)+\text{PD}[m].$$

This choice of the order of $u$ and $w$ is the one that has been implicit
in this paper.\\

We may choose to assign the element
$$\spinc_u(\x)-\frac{1}{2}\text{PD}[m]=\spinc_w(\x)+\frac{1}{2}\text{PD}[m]
\in \frac{1}{2}\text{PD}[m]+\SpinC(Y_0(K)),$$ to the generators
$\x$, which will be for the sake of symmetry as in
\cite{longitude}, or just choose to assign the ``higher'' $\SpinC$
structure $\spinc_u(\x)$ to any $\x \in \Ta \cap \Tb$. Unlike
\cite{longitude} where for  most of the paper we averaged the
$\SpinC$ structures, here we will usually use $\spinc_u(\x)$ as
the $\SpinC$ structure assigned to $\x$. Note that by abuse of
notation,  we are using the notations $\Ta$ and $\Tb$ for the tori
$$\ai_1 \times ...\times \ai_g \times \mu,\text{ and }
l\times \bi_2 \times ...\times \bi_g \times \lambda.$$

Similarly we will obtain the following theorem:

\begin{thm}
There is a natural filtration of the complex $\widehat{CFL}(Y,K)$
respecting the grading by the relative $\SpinC$ structures
$\relspinc \in \RelSpinC(Y,K)=\SpinC(Y_0(K))$ given as
$$\widehat{CFL}_{-}(Y,K,\relspinc)\xra{l_{-}}
\widehat{CFL}_{0}(Y,K,\relspinc)\xra{l_{0}}
\widehat{CFL}_{+}(Y,K,\relspinc),
$$
whose homotopy type is a topological invariant of the knot $(Y,K)$
and the relative $\SpinC$ structure $\relspinc$. Furthermore,
there is an exact sequence
$$0\lra \widehat{CFL}_{0}(Y,K,\relspinc)\xra{l_{0}}
\widehat{CFL}_{+}(Y,K,\relspinc) \xra{\partial_l'}
\widehat{CFL}_{-}(Y,K,\relspinc)\lra 0,
$$
which is also an invariant of the knot $K$. Moreover, if we let
$\partial_l$ be the sum $\partial_l^++\epsilon \partial_l'$, where
$\partial_l^+$ is the boundary map of the complex
$\widehat{CFL}(Y,K)$, then the complex $\widehat{CFL}(Y,K)$,
equipped with $\partial_l$, will give the Floer homology of the
three-manifold $Y_0(K)$. Finally the subcomplex
$\widehat{CFL}_{-}(Y,K)$ has the same chain homotopy type as the
(non-filtered) complex $\widehat{CFK}(Y,K)$.
\end{thm}

\section{Basic properties}

We start with an investigation of the effect of a change in the
orientation of a knot $K$ in the above filtration of
$\widehat{CFK}(Y,K)$. In contrast with the usual hat theory, where
the two points $u$ and $w$ basically played the same role and the
symmetry was relatively easy to show, here we are giving an
essential order to these two marked points, and an essential role
to the orientation of the knot $K$.\\

\begin{figure}
\mbox{\vbox{\epsfbox{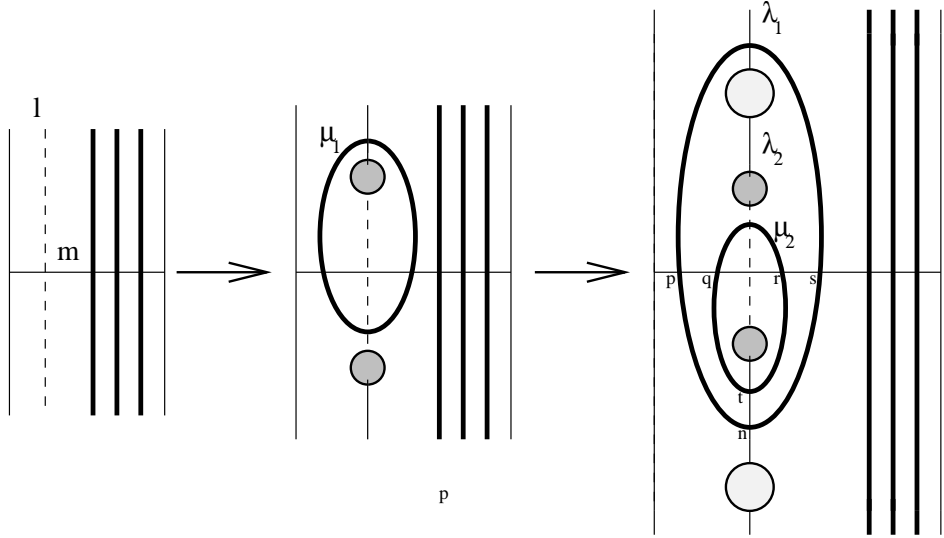}}}
\caption{\label{fig:orientation}
{We may start from a typical Heegaard diagram for the knot $(Y,K)$,
change it to an \emph{allowed} diagram for the pair $(Y,K)$ and then
modify it more to get an \emph{allowed} diagram for $(Y,-K)$.
}}
\end{figure}

Let us start with a standard diagram, not necessarily of the type
discussed in the previous section, for the null-homologous knot
$(Y,K)$. We may run the process of adding a handle in order to get
the appropriate diagram for defining the filtration for $(Y,K)$.
This process will be different if we were going to add the handle
associated with the knot $K$ with the reverse orientation.
However, we may continue the process with this new Heegaard
diagram, as is shown in figure~\ref{fig:orientation}, to add a
second handle, and further modify the longitude of $(Y,-K)$ and
get the
appropriate Heegaard diagram.\\

We will compare the Floer homology of the Heegaard diagram, after
the first step, which gives the filtration of
$\widehat{CFK}(Y,K)$, with that of the final Heegaard diagram,
which gives the filtration on $\widehat{CFK}(Y,-K)$. Let us denote
the $\alpha$ curves added in the first and in the second step by
$\mu_1$ and $\mu_2$ respectively, and similarly denote the new
$\beta$ curves by $\lambda_1$ and $\lambda_2$. Moreover, denote
the \emph{allowed} Heegaard diagram of $(Y,K)$ obtained after the
first step, by $H_1$ and denote the allowed
Heegaard diagram of $(Y,-K)$ obtained after the second step by $H_2$.\\

Some of the special intersection points of these new curves are
named in the diagram of figure~\ref{fig:orientation} by the
letters $n,p,q,r,s$ and $t$. There is a special rectangular region
with vertices $p,q,t$ and $n$ which we denote by $\Delta_L$, and
another rectangular region $\Delta_R$ with the
vertices $r,s,n$ and $t$.\\

Let us partition the set of generators of the complex associated
with  $H_1$ into sets $\mathcal{A},\mathcal{B}$, and
$\mathcal{C}$. Here $\mathcal{A}$ consists of those generators
which are of type A, $\mathcal{B}$ consists of generators of type
B, and $\mathcal{C}$ consists of type C generators.\\

Any generator of type $A$ in $H_2$ will contain the intersection
point $r$. Since we always have to choose one of $t$ and $n$, and
$t$ lies on the same $\alpha$ circle as $r$, the generators of
type $A$ are forced to contain the pair of points $\{r,n\}$. Let
$\mathcal{A}_2$ be set of generators  of $H_2$ of type A. Then the
above discussion shows that there is a one-to-one and surjective
map
$$i_\mathcal{A}:\mathcal{C}\lra \mathcal{A}_2,$$
defined by $i_\mathcal{A}(\{s,\bullet\})=\{r,n,\bullet\}$.
Similarly, if $\mathcal{C}_2$ denotes the set of generators of
$H_2$ of type C, then any of its elements will contain $q$ and
$n$, and there is a correspondence
$$i_\mathcal{C}:\mathcal{A}\lra \mathcal{C}_2,$$
$$i_\mathcal{C}(\{p,\bullet\})=\{q,n,\bullet\}.$$

The generators of type B are more interesting, since they all
share the intersection point $t$. This means that any generator of
$H_1$ may be completed to a generator of $H_2$ by just adding $t$.
As a result we get a correspondence
$$i_\mathcal{B}:\mathcal{A}\cup \mathcal{B} \cup \mathcal{C}
\lra \mathcal{B}_2,$$
$$i_\mathcal{B}(\{\bullet\})=\{t,\bullet\}.$$

If one investigates the effect of the boundary map (differential
map) on the generators of the same type, then the above maps will
further respect the boundary maps, since for example the domain of
a disk connecting $\x$ to $\y$ in $H_1$, is isomorphic to the
domain of the disk connecting $i_\mathcal{B}(\x)$
to $i_\mathcal{B}(\y)$ in $H_2$.\\

Now consider the disks connecting an element $\x$ of
$\mathcal{A}_2$ to an element $\y$ in $\mathcal{B}_2$ or
$\mathcal{C}_2$. Any such disk, if it contributes to the boundary
map, will have non-negative coefficients in all
of the domains.\\

The intersection point $r$ appears in $\x$ and not in $\y$, and
three of the regions which have $r$ as a corner, get a zero
coefficient. As a result the one on the lower right side of $r$ in
figure~\ref{fig:orientation} will get a coefficient equal to $1$.
If $s$ does not appear in $\y$ then a similar comparison of
coefficients around $s$ shows that the region on the upper right
side of $s$ is forced to have a coefficient of $-1$, which is a
contradiction. So $\y$ has to be of type B and in fact the image
of an element of $\mathcal{C}$. Let us assume that
$\x=i_\mathcal{A}(\z)$ and $\y=i_\mathcal{B}(\w)$, with $\z,\w \in
\mathcal{C}$. Then the disk between $\x$ and $\y$ will have a
domain which is the the disjoint union of the rectangle $\Delta_R$
with another domain. An argument similar to the arguments used by
Ozsv\'ath and Szab\'o in proving the connected sum formulas of
\cite{OS-knot} and \cite{OS-3m1,OS-3m2}, may be used to show that
the boundary maps going from $\mathcal{A}_2$ to $\mathcal{B}_2$ is
a perturbation of the map $f_\mathcal{A}=i_\mathcal{B}\circ
i_\mathcal{A}^{-1}$. The total boundary map on $\mathcal{A}_2$
will be chain homotopic to
$$\partial_\mathcal{A}=i_\mathcal{A}\circ \partial \circ i_\mathcal{A}^{-1}
+f_\mathcal{A}=i_\mathcal{A}\circ \partial \circ i_\mathcal{A}^{-1}
+i_\mathcal{B}\circ i_\mathcal{A}^{-1},$$
where $\partial$ denotes the boundary map of $\widehat{CFK}(Y,K)$.\\

Similarly the differentials on type B differ from that of
$\widehat{CFK}(Y,K)$ by a factor
$$f_\mathcal{B}=
i_\mathcal{C}\circ i_\mathcal{B}^{-1}|_{i_\mathcal{B}(\mathcal{A})}.$$

The above discussion may be summarized in the following diagram:
\begin{equation}\label{eq:diagram}
\begin{array}{cc}
&\mathcal{A}\\
&\Big\downarrow\vcenter{\rlap{$d_1$}}\\
&\mathcal{B}\\
&\Big\downarrow\vcenter{\rlap{$d_2$}} \\
&\mathcal{C}
\end{array}
\ \ \ \ \ \ \ \ \ \ \lra \ \ \ \ \ \ \ \ \
\begin{array}{cccccc}
&\mathcal{C}\\
&\Big\downarrow\vcenter{\rlap{$I_\mathcal{C}$}}\\
&\mathcal{C}& \xleftarrow{d_2}&\mathcal{B}& \xleftarrow{d_1}&\mathcal{A}\\
&&&&& \Big\downarrow\vcenter{\rlap{$I_\mathcal{A}$}}\\
&&&&&\mathcal{A}
\end{array}
\end{equation}

The maps $I_\mathcal{A}$ and $I_\mathcal{C}$ are isomorphisms and
$d_1$ and $d_2$ denote just some of the boundary maps. Potentially
there can be maps within each of the sets
$\mathcal{A},\mathcal{B}$ or $\mathcal{C}$, or we can have a
boundary map going from $\mathcal{A}$ to $\mathcal{C}$ on the left
hand side, or in the middle row of the
right hand side. These potential maps are dropped from
diagram~(\ref{eq:diagram}) for simplicity. \\

Since the restriction of the boundary map of $H_1$ to intersection
points of type A (i.e. $\mathcal{A}$) gives a complex isomorphic
to $\widehat{CFK}_{-}(Y,K)$,  We may deduce that up to a possible
total shift in the absolute grading by the Maslov index, there is
a homotopy equivalence
$$\widehat{CFK}_{-}(Y,K)\simeq \widehat{CFK}_{-}(Y,-K).$$

The isomorphism in zero level of the filtration (i.e.
$\widehat{CFK}_0(Y,\pm K)$) is just an algebraic fact about the
two complexes above. And of course we already know the equivalence
of the chain homotopy types of the full complexes
$\widehat{CFK}_{+}(Y,K)$ and $\widehat{CFK}_{+}(Y,-K)$.

When the relative $\SpinC$ structure $\relspinc \in
\RelSpinC(Y,K)$ extends a torsion $\SpinC$ structure $\spinc \in
\SpinC (Y)$, besides the above isomorphism of the homotopy types
of the filtered complexes $\widehat{CFK}(Y,K,\relspinc)$ and
$\widehat{CFK}(Y,-K,\relspinc)$, we may compare the absolute
Maslov gradings of the two complexes.
 Comparing the Maslov indices of the intersection points with the
definition of Ozsv\'ath and Szab\'o of the absolute Maslov grading
in \cite{OS-knot}, we derive the following:

\begin{prop}
If $(Y,K)$ is a null-homologous knot and if $(Y,-K)$ denote the
same knot with the opposite orientation, then for any relative
$\SpinC$ structure $\relspinc \in \RelSpinC(Y,K)$, the filtered
complexes  $\widehat{CFK}(Y,K,\relspinc)$ and
$\widehat{CFK}(Y,-K,\relspinc)$ have the same homotopy type.
Moreover if $\relspinc$ extends a torsion $\SpinC$ structure
$\spinc \in \SpinC(Y)$, then the homotopy equivalence shifts the
absolute Maslov grading by $-2k$, where
$$k=\frac{1}{2}\langle c_1(\relspinc),[\widehat{F}]\rangle,$$
for a capping $\widehat{F}$ of a Seifert surface $F$ for the oriented knot
$K$.  In particular
$$H_d(\widehat{CFK}_{\bullet}(Y,K,\relspinc))\simeq
H_{d-2m}(\widehat{CFK}_{\bullet}(Y,-K,\relspinc)),
$$
where $\bullet \in \{-,0,+\}$.
\end{prop}

\begin{proof}
Note that the whole complex $\CFK_+(Y,K)$ is precisely the complex
$\CFK(Y,K)$ in \cite{OS-knot} and the absolute Maslov grading of
generators is induced using this identification. Using this piece
of information and the corresponding isomorphism
$$\CFK(Y,K)\simeq \CFK(Y,-K)$$ in proposition 3.8 of
\cite{OS-knot}, the proof is just to verify the statement via an
algebraic comparison.
\end{proof}

The next thing we want to consider in this section is the effect
of the maps in the long exact sequences on the above filtration.

We remind the reader of the basic setup of the long exact
sequences in \cite{OS-knot} and \cite{OS-3m2}. Recall that for any
framed knot $\gamma$ in the three manifold $Y$, and for any fixed
$\SpinC$ structure $\spinct \in \SpinC(W_\gamma(Y))$ on the
four-manifold $W=W_\gamma(Y)$, counting holomorphic triangles
gives a map
$$\widehat{f}_{W,\spinct}:\widehat{CF}(Y) \lra \widehat{CF}(Y_\gamma),$$
where $Y_\gamma$ is obtained from $Y$ by a surgery on the framed
knot $\gamma$. For the definition of the $4$-manifold
$W_\gamma(Y)$ we refer the reader to \cite{OS-4m}. When we have a
knot $(Y,K)$ whose linking number with the framed knot $\gamma$ is
zero, we may choose a Seifert surface $F$ for $K$ which is
disjoint from $\gamma$. Correspondingly, we will get a filtration
induced
by $F$ on both $\widehat{CF}(Y)$ and $\widehat{CF}(Y_\gamma)$.\\

Ozsv\'ath and Szab\'o showed in \cite{OS-knot} that in fact the
map $\widehat{f}_{W,\spinct}$ respects this  filtration. As a
result, from the long exact sequence for the three-manifolds, they
obtain the following exact sequence for knots:
$$...\ra \widehat{HFK}(Y_{-1}(\gamma),K,k) \xra{f^1}
 \widehat{HFK}(Y_{0}(\gamma),K,k) \xra{f^2}
 \widehat{HFK}(Y,K,k) \xra{f^3} ...,$$
for any integer $k\in \mathbb{Z}$, where we use their notation:
$$\widehat{HFK}(Y,K,k)=\bigoplus_{\substack{
\relspinc \in \RelSpinC(Y,K)\\  \langle
c_1(\relspinc),[\widehat{F}]\rangle =2k}}
\widehat{HFK}(Y,K,\relspinc).$$

We claim that the following is also true:

\begin{thm}
Let $(Y,K)$ be a knot and let $\gamma$ be a framed knot in $Y$
which is disjoint from $K$ and has zero linking number with it.
Then for a correct choice of the Seifert surface, the counts of
holomorphic triangles above respects the $\{-,0,+\}$ filtration,
and for each integer $k\in \mathbb{Z}$ we obtain the exact
sequences:
\begin{equation}
\begin{split}
...\ra H_*(\widehat{CFK}_{\bullet}(Y_{-1}(\gamma),K,k))
&\xra{f^1_{\bullet}}
 H_*(\widehat{CFK}_{\bullet}(Y_{0}(\gamma),K,k)) \xra{f^2_{\bullet}}\\
&\ra H_*(\widehat{CFK}_{\bullet}(Y,K,k)) \xra{f^3_{\bullet}}... ,
\end{split}
\end{equation}
for any of the filtration levels $\bullet \in \{-,0,+\}$.  The
maps $f^{\bullet}_1$ and $f^{\bullet}_2$, when the groups are
graded, each will lower the absolute grading by $\frac{1}{2}$, and
$ f^{\bullet}_3$ will not increase the absolute grading.
\end{thm}

\begin{proof}
The proof is quite straight forward. In the triple Heegaard
diagrams used to connect the  knot $(Y,K)$ to the knot
$(Y_\gamma,K)$, the argument of the second section shows, by
examining the positivity of the coefficients of potential
holomorphic disks in the domains around the two marked points,
that the image of a generator  of type C, can only have type C,
and that the image of a generator of type B is either of type B,
or of type C. As a result the filtration of $\widehat{CFK}(Y,K,k)$
as
$$\widehat{CFK}_{-}(Y,K,k)\subset \widehat{CFK}_{0}(Y,K,k)\subset
\widehat{CFK}_{+}(Y,K,k)$$ is preserved under the connecting
homomorphisms. This completes the proof.
\end{proof}

\section{Gluing along the knots I}

Suppose that $(Y_1,K_1)$ and $(Y_2,K_2)$ are two oriented knots.
Then there is a unique framing on the boundary of each of the
three-manifolds with boundary $Y_1\setminus \text{nd}(K_1)$ and
$Y_2\setminus \text{nd}(K_2)$, determined by a meridian $m_i$ of
each of the knots $K_i$, and a longitude $l_i$, such that each
$l_i$ corresponds to a zero surgery on $K_i$,
$i=1,2$ (i.e. $l_i$ has zero linking number with $K_i$).\\

We may glue these two three-manifolds with boundary, by
identifying $l_1$ with $l_2$ and $m_1$ with $-m_2$ (as oriented
curves), along their torus boundaries. In fact, the result will be
$$(Y_1\setminus \text{nd}(K_1))\cup_T(Y_2\setminus \text{nd}(K_2)),$$
where $T$ is the identified  boundary of $\text{nd}(K_1)$ and
$\text{nd}(K_2)$. This will give us a new three-manifold $Y$ and
inside this three-manifold, the knots $K_1$ and $K_2$ will be
identified to give a knot $(Y,K)$. Sure enough, $K$ is
null-homologous in $Y$. We will denote this operation by writing
$$(Y,K)=(Y_1,K_1)\|(Y_2,K_2).$$

In this section we will find out how the (non-filtered) Floer homology of
$(Y,K)$ is related to filtered Floer homologies associated with
$(Y_i,K_i)$, $i=1,2$. In fact, we will prove
the following:\\

\begin{thm}\label{thm:gluing}
Suppose that $(Y_1,K_1)$ and $(Y_2,K_2)$ are two null-homologous
knots and that $(Y,K)=(Y_1,K_1)\|(Y_2,K_2)$ is obtained by gluing
$(Y_1,K_1)$ to $(Y_2,K_2)$ as above. Then for any relative
$\SpinC$ structure
$$\relspinc_1\#\relspinc_2 \in \RelSpinC(Y,K)=\RelSpinC(Y_1,K_1)\oplus
\RelSpinC(Y_2,K_2),$$ the Heegaard Floer homology of $(Y,K)$ in
the $\SpinC$ structure $\relspinc_1\#\relspinc_2$ will be given by
the quotient complex

\begin{displaymath}
\sC=\frac{\Big [ \frac{\CFL(Y_1,K_1,\relspinc_1)\otimes
\CFL(Y_2,K_2,\relspinc_2)} {\CFL_0(Y_1,K_1,\relspinc_1)\otimes
\CFL_0(Y_2,K_2,\relspinc_2)}\Big ] }{\Large{\sim}}.
\end{displaymath}
Here the equivalence relation $\sim$ is induced by the isomorphism
of the subcomplexes
\begin{displaymath}
\begin{split}
\rho: \CFL_-(Y_1,K_1,\relspinc_1)&\otimes \frac{\CFL(Y_2,K_2,\relspinc_2)}
{\CFL_0(Y_2,K_2,\relspinc_2)}\\
&\lra
\frac{\CFL(Y_1,K_1,\relspinc_1)}{\CFL_0(Y_1,K_1,\relspinc_1)}\otimes
\CFL_-(Y_2,K_2,\relspinc_2),
\end{split}
\end{displaymath}
where $\rho$ is induced by the isomorphisms $\CFL_-\simeq \frac{\CFL}{\CFL_0}$.
\end{thm}
The rest of this section is devoted to a proof of this theorem. \\

\begin{proof}
Suppose that $H_i$ is a Heegaard diagram consisting of a surface
$\Sig_i$ of genus $g_i$ and two sets of curves
$$\alphas^i=\{\ai^i_1,...,\ai^i_{g_i}\}, \text{ and}$$
$$\betas^i_0=\{\bi^i_2,...,\bi^i_{g_i}\},$$
together with two special curves $m_i$ and $l_i$ representing the
meridian and the longitude of $K_i$, for $i=1,2$. Moreover assume
that $l_i$ and
$m_i$ meet each other transversely in a single point $z_i$.\\

In order to obtain a Heegaard diagram for $(Y,K)$, connect the surfaces
$\Sig_1$ and $\Sig_2$ by a handle which connects $z_1$ to $z_2$. Connect
$l_1$ to $l_2$ using this handle so that it represents the trivial knot
$K_1 -K_2$. Also, connect
$m_1$ to $m_2$ so that the resulting closed curve represents the trivial
curve $m_1+m_2$. Then the resulting surface $\Sig$ together with the
$g_1+g_2$ tuples of curves:
$$\alphas=\alphas^1 \cup \alphas^2, \text{ and}$$
$$\betas=\betas_0\cup\{m_1\#m_2\}=\betas^1_0 \cup \betas^2_0\cup \{l_1\#l_2\}
\cup \{m_1 \# m_2\}$$
will form a Heegaard diagram for $(Y,K)$ with the meridian $m=m_1\#m_2$.
Denote this Heegaard diagram by $H_1*H_2$.\\

\begin{figure}
\mbox{\vbox{\epsfbox{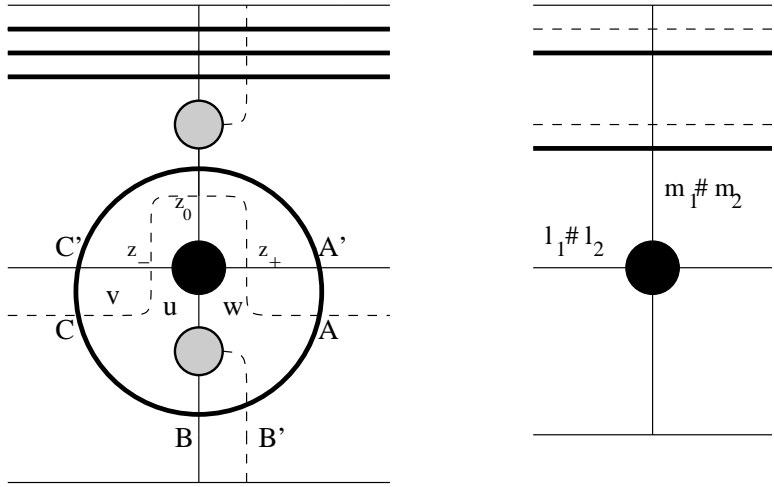}}} \caption{\label{fig:surgery}
{An allowed Heegaard diagrams (with respect to the longitude
theory $\widehat{CFL}$) for $(Y_1,K_1)$ (on the left) and a weakly
admissible Heegaard diagram for $(Y_2,K_2)$ (on the right) may be
connected by a handle, denoted by the black circle, to obtain a
Heegaard diagram for $(Y,K)$. }}
\end{figure}

We may assume that $H_1$ and $H_2$ are any specific Heegaard
diagrams for the the knots $(Y_1,K_1)$ and $(Y_2,K_2)$,
respectively. In particular, we may assume that $H_1$ is already
modified, so that it is a suitable (allowed) Heegaard diagrams for
defining the filtered version of $\widehat{CFL}(Y_1,K_1)$. We may
also assume that $H_2$ is a weakly admissible diagram, which may
be used for a computation of $\widehat{CFL}(Y_2,K_2)$. It is then
implied that $H=H_1*H_2$ is weakly admissible. The Heegaard
diagram $H$
 may still not be an allowed Heegaard diagram for computing $\CFK
(Y,K)$. However, this Heegaard diagram may be used to
compute the non-filtered version of $\CFK(Y,K)$.\\

The final Heegaard diagram $H$ which defines the complex $\widehat{CFK}(Y,K)$
will look like figure~\ref{fig:surgery} around the attached handle. In this figure,
the attaching circles of the connecting handle of $\Sig_1$ and $\Sig_2$ are
 denoted by the black circles
on the two sides, and each other pair of circles of the same color represents another
handle. Moreover, the regular curves denote the elements of $\betas$, while the
bold curves represent the elements of $\alphas$.\\

One may consider another set of $g_1+g_2$ curves, denoted by
$\thetas$, on the surface $\Sigma$. The set $\thetas$ will
formally be the union
$$\thetas=\betas^1_0\cup \alphas^2 \cup \{\mu\}
=\deltas_0\cup\gammas \cup \{\mu\}$$
$$=\{\delta_2,...,\delta_{g_1}\}\cup\{\gamma_1,...,\gamma_{g_2}\}\cup\{\mu\},$$
where $\gammas$ will consist of  Hamiltonian isotopes of the
curves in $\alphas^2$ so that the intersections of all pairs of
curves are transverse, and consisting of cancelling intersections.
Similarly, we assume that the curves in $\deltas_0$ are
Hamiltonian isotopes of those in $\betas^1_0$. The special curve
$\mu$ will be  a Hamiltonian isotope of the longitude $l_1$ on
$\Sig_1$ side of the picture, such that it cuts the curve
$l_1\#l_2$ in a pair of cancelling intersection points $z_+$ and
$z_-$, and cuts the meridian $m_1\#m_2$ in a single point $z_0$,
and such that all the intersection points are very close to the
attaching circle of the handle connecting $\Sig_1$ to $\Sig_2$. We
may choose all the above curves so that the picture is as
illustrated in figure~\ref{fig:surgery}, with the dotted curves
representing the curves in  $\thetas$. We choose $\mu$ so that
it stays very close to $l_1$ after leaving $z_+$ for $z_-$. Then
from $z_-$ to $z_0$ and from $z_0$ to $z_+$ it stays very close
to the $\alpha$-curve of figure~\ref{fig:surgery}, which we assume
to be $\alpha_1$. In particular we assume that $z_-$ and $z_+$ are
very close to the intersection points in $\alpha_1\cap (l_1\#l_2)$
and that $z_0$ is very close to $\alpha_1\cap (m_1\#m_2)$.\\

Correspondingly, we may consider the triple Heegaard diagram
$$\overline{H}=(\Sig,\alphas,\thetas,\betas;u,w,v),$$
consisting of the above three sets of  curves on $\Sig$, together
with three marked points $u,w$ and $v$.  The triple Heegaard
diagram $\Hbar$ gives a map $\mathcal{F}$, through a count of
holomorphic triangles associated with this Heegaard diagram
$\overline{H}$. In order to obtain $\sF$, we will count only those
triangles such that their domain has coefficient zero at $u$ and
$w$, and  $v$. The map $\sF$ goes from the complex
$$\widehat{CF}(\alphas,\thetas)\otimes \widehat{CF}(\thetas,\betas)$$
to the complex $$\widehat{CF}(\alphas,\betas)=\widehat{CFK}(Y,K).$$
Here, $\widehat{CF}(\alphas,\thetas)$ denotes the chain complex associated with
the Heegaard diagram
$$(\Sig,\alphas,\thetas; u=w, v),$$
where  the "$=$" signs between $u$ and $w$ denotes that they are in the same domain
of the Heegaard diagram, when we get rid of $\betas$. \\

Similarly, for $\widehat{CF}(\thetas,\betas)$ we will use the Heegaard diagram
$$(\Sig,\thetas,\betas;u=w,v).$$
The "equality " $w=u$ is less trivial in this case, compared to
the previous equality. Note the difference between the situation
here compared to that of the holomorphic triangle map which
appeared in the proof of connected sum formulas
(\cite{OS-3m2,OS-knot}). \\

For any pair of relative
$\SpinC$ structures $\relspinc_1 \in \RelSpinC(Y_1,K_1)$ and
$\relspinc_2 \in \RelSpinC(Y_2,K_2)$, there is a unique relative $\SpinC$ structure
$$\relspinc=\relspinc_1 \# \relspinc_2 \in \RelSpinC(Y,K)$$
which extends $\relspinc_1$ and $\relspinc_2$. This should be  compared with the
$\mathbb{Z}$-many of pairs of relative $\SpinC$ structures resulting in the same total
relative $\SpinC$ structure, when forming the connected sum of the two knots
$(Y_1,K_1)$ and $(Y_2,K_2)$.\\

In order to understand the chain complex
$\widehat{CFK}(Y,K)=\widehat{CF}(\alphas,\betas)$, we will perform
a study of the two chain complexes $\widehat{CF}(\alphas,\thetas)$
and $\widehat{CF}(\thetas,\betas)$, followed by analyzing the
effect of the map
$\mathcal{F}$. These pieces of information will be enough to prove the above theorem.\\

First of all note that the Heegaard diagram $H^1=(\Sig,\alphas,\thetas,w,v)$ represents
the connected sum $V_1$ of the three-manifold $(Y_1)_0(K_1)$,
 obtained from $Y_1$ by a zero surgery
on the knot $K_1$, with the three manifold $\#^{g_2}S^1\times S^2$.
For any relative $\SpinC$ structure
$\relspinc_1\in \RelSpinC(Y_1,K_1)$ we may consider the $\SpinC$ structure
$\relspinc_1 \# \relspinc_0$ on $V_1$, where $\relspinc_0$ is the $\SpinC$ structure on
$\#^{g_2}S^1\times S^2$ with the property that $c_1(\relspinc_0)$ represents the trivial
homology class. Similar to the methods used in \cite{OS-3m2,OS-knot}, we may couple
any generator $\x$ of the Heegaard diagram
$$\overline{H}_1=(\Sig_1,\alphas^1,\betas^1_0\cup\{\mu\},w,v),$$
corresponding to the relative $\SpinC$ structure $\relspinc_1$
(i.e. $\spinc_w(\x)=\relspinc_1$), with the generator
$\Theta_2$ of the complex associated with the Heegaard diagram
$(\Sig_2,\alphas^2,\gammas,u=w)$ corresponding to the
top generator of the homology
$$\widehat{HF}(\alphas^2,\gammas)=
\widehat{HF}(\#^{g_2}S^1\times S^2,\relspinc_0),$$ to obtain a
corresponding generator of the Heegaard diagram $H^1$ in $\SpinC$
class $\relspinc_1\#\relspinc_2$.\\

 As a result, the chain complex $C_1$ associated with the Heegaard diagram $\Hbar_1$
may be embedded via the above map into the chain complex $\CF(H^1)$:
$$p_1:C_1(\relspinc_1) \xra{\times \Theta_2} \CF(H^1,\relspinc_1\#\relspinc_0).$$
We are introducing the notation $\sC(\relspinc)$ for the part of a
chain complex $\sC$ corresponding to the $\SpinC$ class
$\relspinc$.\\

The complex $C_1$ may be identified with the filtered chain
complex $\CFL_\bullet(Y_1,K_1)$, and under this identification
$C_1(\relspinc_1)=\CFL_+(Y_1,K_1,\relspinc_1)$.

On the other side, the story is slightly more complicated. Here, in the Heegaard diagram
$H^2=(\Sig,\thetas,\betas,u=w,v)$, we may  consider the two sets of curves
$\deltas_0$ and $\betas^1_0$, which are Hamiltonian isotopes of one another. As a result,
there will be a map to the chain complex associated with $H^2$, from the chain
complex associated with the Heegaard diagram
$$\Hbar_2=
(\Sig_2',\{\mu\}\cup \gammas, \betas^2_0 \cup \{l_1\#l_2\} \cup \{m_1\#m_2\},
u=w,v).$$
Here $\Sig_2'$ is obtained from $\Sig_2$ by making a connected sum with a torus $T$, in
the same way that the surface $\Sig$ is obtained from $\Sig_2$ by taking a connected sum
of it with $\Sig_1$. The last two curves $l_1\#l_2$ and $m_1\#m_2$ are in fact the
curves induced from $\Sig_1\#\Sig_2$ to $T\#\Sig_2$. The above map
will be given by taking a product with the top generator of the Floer homology of
$\#^{g_1-1}S^1\times S^2$ obtained from the $(g_1-1)$ tuples of curves $\deltas_0$
and $\betas^1_0$. We will denote the chain complex associated with $\Hbar_2$ by $C_2$.\\

If we drop the marked points $u,w$ and $v$ from the Heegaard
diagram $\Hbar_2$, we obtain a Heegaard diagram for
$(Y_2)_0(K_2)$, and the maps $\spinc_u$ and $\spinc_w$ from the
generators of $C_2$ to $\RelSpinC(Y_2,K_2)$ will differ from each
other by a factor $\PD[m_2]$. In fact for our choice of $u$ and
$w$ we will have $\spinc_u(\x)=\spinc_w(\x)+\PD[m_2]$. We may
assign a $\SpinC$ structure $\relspinc_2$ to any such generator,
using the map $\spinc_u$, and decompose $C_2$ into summands
$C_2(\relspinc_2)$ accordingly.\\

We may re-draw the local configuration of the curves in the Heegaard diagram $\Hbar_2$,
which results in the picture shown in figure~\ref{fig:diagramH3}. The above discussion
shows that there is a second embedding:
$$p_2:C_2(\relspinc_2)\xra{\Theta_1 \times } \CF(\Hbar_2, \relspinc_0\# \relspinc_2).$$
We are of course abusing the notation by denoting the special
$\SpinC$ structure of $\#^{g_1-1}S^1 \times S^2$ by the same
symbol $\relspinc_0$, which was used for
the special $\SpinC$ structure on $\#^{g_2}S^1\times S^2$.\\

Composing these embeddings with the count $\mathcal{F}$ of the holomorphic triangles,
we obtain a chain map, still denoted by $\mathcal{F}$:
$$\mathcal{F}:(C_1,\relspinc_1) \otimes (C_2,\relspinc_2) \lra
\widehat{CFK}(Y,K,\relspinc_1\# \relspinc_2),$$
$$\relspinc_1 \in \RelSpinC(Y_1,K_1), \relspinc_2 \in \RelSpinC(Y_2,K_2). $$

Any generator of the complex $C_2$ will contain an intersection point
on the image of the special curve $\mu$. There are three intersection points between
$\mu$ and the two curves $\ell=l_1\#l_2$ and $m=m_1\#m_2$. The ones on $\ell$ are
denoted by $z_+$ and $z_-$, and the intersection with $m$ is denoted by $z_0$.
We choose them so that the trivial disk between the generators $\{z_+,\bullet\}$
and $\{z_-,\bullet\}$ goes from the first one to the second.
Correspondingly, we may divide the generators of $\widehat{CF}(\Hbar_2)$ into the
 groups $\zZ_+,\zZ_-$ and $\zZ_0$. \\

It is clear that the complex $C_2$ may be identified with the filtered complex
$\CFL(Y_2,K_2)$ and that there is a natural bijection
$$j:\CF(\Hbar_2) \lra \CFL(Y_2,K_2),$$
which reduces to maps
$$j_-:\sZ_- \lra \CFL_-(Y_2,K_2),$$
$$j_0:\sZ_-\cup\sZ_0 \lra \CFL_0(Y_2,K_2) \ \ \ \ \text{ and}$$
$$j_+:\sZ_-\cup \sZ_0 \cup \sZ_+ \lra \CFL_+(Y_2,K_2).$$

\begin{figure}
\mbox{\vbox{\epsfbox{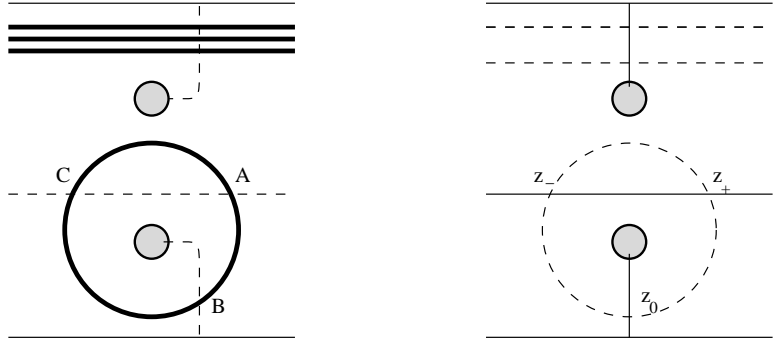}}}
\caption{\label{fig:diagramH3} { Associated with the pair
$(\alphas,\thetas)$, we will get a
  connected  sum of the diagram on the left with a Heegaard diagram for
$\#^{g_2}S^1\times  S^2$, and associated with the pair
  $(\thetas,\betas)$ we will get a connected sum of a Heegaard diagram for
  $\#^{g_1-1}S^1\times S^2$ with the Heegaard diagram on the right.
}}
\end{figure}

Under the above identifications
$$C_1\simeq \CFL(Y_1,K_1),\text{ and }$$
$$C_2\simeq \CFL(Y_2,K_2),$$
suppose that $\x \otimes \y \in C_1 \otimes C_2$ belongs to
$\CFL_0(Y_1,K_1) \otimes \CFL_0(Y_2,K_2)$. This means that the generator
$\y$  contains one of the intersection points $z_-$ or $z_0$, and that the intersection
point $A$ of figure~\ref{fig:surgery} is not included in the generator $\x$.
It is not hard to check, by examining the possible coefficients of the formal domain
associated with any holomorphic  triangle, that the image $\sF(\x \otimes \y)$ is
trivial. In fact, if $\phi$ is a triangle connecting $\x \times \Theta_2$,
$\Theta_1 \times \y$ and $\w$, with $\w \in \CFK(Y,K)=\CF(H)$, then $\phi$ can not have
non-negative coefficients in all the domains appearing in figure~\ref{fig:surgery},
and zero  coefficients at $u,w$ and $v$.\\

As a result of the above discussion we see that $\CFL_0(Y_1,K_1) \otimes \CFL_0(Y_2,K_2)$
is included in the kernel of the map $\sF$. So we may think of $\sF$ as a map
$$\sF:\frac{\CFL(Y_1,K_1)\otimes \CFL(Y_2,K_2)}{\CFL_0(Y_1,K_1)\otimes
\CFL_0(Y_2,K_2)} \lra \CFK(Y,K).$$\\

Now suppose that $\x_- \otimes \y_+$ is a generator in
$$\CFL_-(Y_1,K_1) \otimes \frac{\CFL(Y_2,K_2)}{\CFL_0(Y_2,K_2)}$$
which corresponds to a coupling of a generator $\x_-$ of type C
and a generator $\y_+$ of type A. Denote the type A generator of
the complex $C_1$ which corresponds  to $\x_-$, by $\x_+$.
Similarly, the type C generator of $C_2$ which is associated with
$\y_+$ will be denoted by $\y_-$. We will obtain a generator
$$\x_+ \otimes \y_- \in \frac{\CFL(Y_1,K_1)}{\CFL_0(Y_1,K_1)}\otimes \CFL_-(Y_2,K_2).$$

There is a unique generator $\w \in \CF(H)\simeq \CFK(Y,K)$ with
the property that there is triangle connecting $\x_- \times
\Theta_2$, $\Theta_1 \times \y_+$ and $\w$, whose domain consists
of $g_1+g_2$ disjoint small triangles with a  total area that may
be assumed to be arbitrarily small, if we choose the curves in
$\thetas$ sufficiently close to the corresponding curves in
$\alphas$ and $\betas$. This triangle supports a unique
holomorphic representative. Moreover, the same generator $\w$ has
the extra property that there is another unique holomorphic
triangle $\psi$ with small area (energy) that connects $\x_+
\times \Theta_2$, $\Theta_1 \times \y_-$ and $\w$. As a result
$$\sF(\x_- \otimes \y_+)=\w + \text{ significantly lower energy terms, }$$
and
$$\sF(\x_+ \otimes \y_-)=\w + \text{ significantly lower energy terms. }$$
Using the arguments of Ozsv\'ath and Szab\'o on the energy
filtration (\cite{OS-3m2}~, section 6), we may use this
correspondence to modify the isomorphism
$$\rho:\CFL_-(Y_1,K_1)\otimes \frac{\CFL(Y_2,K_2)}{\CFL_0(Y_2,K_2)}\lra
\frac{\CFL(Y_1,K_1)}{\CFL_0(Y_1,K_1)} \otimes \CFL_-(Y_2,K_2)$$ by
a chain homotopy, so that for any
$$\z=\x_-\otimes \y_+ \in \CFL_-(Y_1,K_1) \otimes \frac{\CFL(Y_2,K_2)}{\CFL_0(Y_2,K_2)},$$
we have $\sF(\z -\rho (\z))=0$. \\

With this observation in hand, it is now enough to show that the
map $\sF$, induced on the quotient complex
$$\sC=\big [
\frac{\CFL(Y_1,K_1)\otimes \CFL(Y_2,K_2)}{\CFL_0(Y_1,K_1)\otimes
\CFL_0(Y_2,K_2)}\big ]/\sim
$$ to $\CFK(Y,K)$ is a bijection, where $\sim$ is induced by the
modified isomorphism $\rho$.\\

But this is quite easy. The map $\sF$ is of the form
$$\sF=\sG+ \text{ significantly lower energy terms, }$$
where the map $\sG$ is defined from
\begin{displaymath}
\begin{split}
\big [ &\frac{\CFL(Y_1,K_1)\otimes
\CFL(Y_2,K_2)}{\CFL_0(Y_1,K_1)\otimes \CFL_0(Y_2,K_2)}\big ]/
\sim \\
&\ \ \  \Big [\CFL(Y_1,K_1) \otimes
\frac{\CFL(Y_2,K_2)}{\CFL_0(Y_2,K_2)}\Big ] \bigoplus \Big [
\frac{\CFL(Y_1,K_1)}{\CFL_0(Y_1,K_1)} \otimes
 \frac{\CFL_0(Y_2,K_2)}{\CFL_-(Y_2,K_2)} \Big ]
\end{split}
\end{displaymath}
to $\CFK(Y,K)$ as follows. The first summand corresponds to a
coupling of an arbitrary element $\x$ of $\CFL(Y_1,K_1)$, with an
element $\y$ of $\sZ_+$. There is a unique generator $\w=\sG(\x
\otimes \y)$ of $\CFK(Y,K)$ such that the three generators $\x
\otimes \Theta_2, \Theta_1 \times \y$ and $\w$ are connected by a
unique triangle $\phi$, whose domain $\sD(\phi)$ consists of small
triangles of arbitrarily small area,
and such that $\phi$ supports a holomorphic disk. \\

The definition of $\sG$ on the second summand is done in a quite
similar fashion. The argument which finds the multiplicities of
the connecting disk is slightly more complicated. An examination
of possible multiplicities for the domain of such a triangle, in
the part of surface that appears in figure~\ref{fig:surgery}, and
just another game of trying to keep the multiplicities positive,
will find two candidates for $\w= \sG(\x \otimes \y)$, where
$$\x \otimes \y \in \frac{CFL(Y_1,K_1)}{\CFL_0(Y_1,K_1)} \otimes
 \frac{\CFL_0(Y_2,K_2)}{\CFL_-(Y_2,K_2)}.$$
However, only one of them has the property that the triangle associated with it has an
arbitrarily small area, once $\mu$ is chosen as it was described before.\\

It is not hard to check that the map $\sG$ is surjective. From the construction it
is clear that our claim is satisfied, i.e.
 $$\sF=\sG+ \text{ significantly lower energy terms. }$$

Now the argument of Ozsv\'ath and Szab\'o in \cite{OS-3m2,
OS-knot} may be followed to conclude that $\sF$ is surjective and
injective from $\sC$ to $\CFK(Y,K)$. This completes the proof of
the theorem.
\end{proof}

\section{Gluing along the knots II}

We may consider several other gluing constructions along the
knots, as well as the operation of taking the connected sum of
$(Y_1,K_1)$ with $(Y_2,K_2)$, and look into the complexes
$\CFL(Y,K)$ and $\CFK(Y,K)$ of the resulting knot. Now that we
have
treated the above case carefully, it is not hard to state an prove some similar results.\\

The first of these results, is the effect of connected sum on $\CFL$. We already know
that (see \cite{OS-knot}) if $(Y_1,K_1)$ and $(Y_2,K_2)$ are null homologous knots
and if $(Y,K)$ is obtained by taking the connected sum of them, then for any
relative $\SpinC$ structure $\relspinc \in \RelSpinC(Y,K)$ there are $\Z$-many
of pairs
$$\relspinc_1\# \relspinc_2\in \RelSpinC(Y_1,K_1)\oplus \RelSpinC(Y_2,K_2)$$
which extend $\relspinc$. The complex $\CFK(Y,K)$ may be described as
\begin{displaymath}
\CFK(Y,K,\relspinc)=\bigoplus_{\substack{\relspinc_i \in \RelSpinC(Y_i,K_i)\\
\relspinc_1\#\relspinc_2 \text{ extends } \relspinc}}\CFK(Y_1,K_1,\relspinc_1)
\otimes \CFK(Y_2,K_2,\relspinc_2).
\end{displaymath}

We claim that $\CFL(Y,K)$ is given by the following theorem.

\begin{thm}
If $(Y,K)$ is obtained as the connected sum of the two knots $(Y_1,K_1)$ and
$(Y_2,K_2)$, and if $\relspinc \in \RelSpinC(Y,K)$ is a relative $\SpinC$ structure
on $(Y,K)$, then
\begin{displaymath}
\begin{split}
&\CFL(Y,K,\relspinc)=\frac{\Big [\bigoplus_{\substack{\relspinc_i \in \RelSpinC(Y_i,K_i)\\
\relspinc_1\#\relspinc_2 \text{ extends } \relspinc}} \sC(\relspinc_1,\relspinc_2) \Big ]}
{\Large \sim}\\
&\sC(\relspinc_1,\relspinc_2)=\frac{
\CFK(Y_1,K_1,\relspinc_1) \otimes \CFK(Y_2,K_2,\relspinc_2)}{\CFK_0(Y_1,K_1,\relspinc_1)
\otimes \CFK_0(Y_2,K_2,\relspinc_2)}\\
&\CFK_-(Y_1,K_1,\relspinc_1)\otimes \frac{\CFK(Y_2,K_2,\relspinc_2)}{\CFK_0(Y_2,K_2,
\relspinc_2)} {\Large{\sim}}\\
&\ \ \ \ \ \ \ \ \ \ \ \ \ \frac{\CFK(Y_1,K_1,\relspinc_1 +
\text{PD}[m_1])}{\CFK_0(Y_1,K_1,\relspinc_1+\text {PD}[m_1])}
\otimes \CFK_-(Y_2,K_2,\relspinc_2 -\text{PD}[m_2]).
\end{split}
\end{displaymath}
Here the isomorphisms inducing $\sim$ are coming from the
connecting isomorphisms $\partial':\frac{\CFK}{\CFK_0}\ra \CFK_-$.
\end{thm}

\begin{proof}
The proof is almost identical to the proof of the above theorem.
In fact, if we change the role of the meridians $m_1$ and $m_2$
with the longitudes $l_1$ and $l_2$, then the Heegaard diagram
$H=H_1 * H_2$ which was used in the proof of
theorem~\ref{thm:gluing} will compute the longitude Floer homology
$\CFL(Y,K)$, where $(Y,K)=(Y_1,K_1)\#(Y_2,K_2)$. The allowed
Heegaard diagrams $\overline{H}_1$ and $\overline{H}_2$, and also
$\Hbar_2$ will compute the filtered chain complexes
$\CFK(Y_1,K_1)$ and $\CFK(Y_2,K_2)$. The equivalence relation
$\sim$ will come from, instead of the isomorphism between the type
A and type C generators of $\CFL$, the isomorphism between type A
and type C generators of $\CFK(Y_i,K_i)$. This last isomorphism
changes the $\SpinC$ structure within each of the complexes
$\CFK(Y_i,K_i)$. That is the reason for the appearance of the
shifts by $+\text{PD}[m_1]$ and $-\text{PD}[m_2]$, which does not
change the total relative $\SpinC$ structure $\relspinc \in
\RelSpinC(Y,K)$, since the total change is by the Poincar\'e dual
of $[m_1-m_2]$.
\end{proof}

Since the effect of a connected sum on the longitude theory is so
closely related to the effect of the gluing along the knots on the
standard complex $\CFK$, we may expect that the converse is also
true. In fact, if we drop our filtration of the complex $\CFL$,
introduced in this paper, we get the following proposition.

\begin{prop}
Suppose that the knot $(Y,K)$ is obtained from $(Y_1,K_1)$ and
$(Y_2,K_2)$, by gluing the knot complements $Y_1 \setminus
\text{nd}(K_1)$ and  $Y_2 \setminus \text{nd}(K_2)$ along their
boundary and identifying the longitudes (resp. meridians) of $K_1$
and $K_2$ (i.e. we have $(Y,K)=(Y_1,K_1)\|(Y_2,K_2)$). Then for
any
$$\relspinc=\relspinc_1\#\relspinc_2 \in
\RelSpinC(Y,K)=\RelSpinC(Y_1,K_1)\oplus \RelSpinC(Y_2,K_2)$$
we will have
$$\CFL(Y,K,\relspinc)=\CFL(Y_1,K_1,\relspinc_1)\otimes \CFL(Y_2,K_2,\relspinc_2).$$
\end{prop}
\begin{proof}
Again the proof is almost identical to the proof of the connected
sum formula for $\CFK(Y,K)$ in \cite{OS-knot}. In the Heegaard
diagram $H=H_1*H_2$ of theorem~\ref{thm:gluing}, $m=m_1\#m_2$ will
be the meridian of the knot $K$, and a parallel copy $l_1'$ of
$l_1$ may be used as the longitude of $K$. Thus, a Heegaard
diagram for computing $\CFL(Y,K)$ may be obtained by trading $m$
for $l_1'$ in the Heegaard diagram $H$. Up to a change in the
names, this is similar to the Heegaard diagram used in section 7
of \cite{OS-knot}. The argument of Ozsv\'ath and Szab\'o may be
followed word by word.
\end{proof}

The next thing we may derive without any extra effort is a gluing
formula for the following construction. Suppose that $(Y_1,K_1)$
and $(Y_2,K_2)$ are as before. We may glue $Y_1 \setminus
\text{nd}(K_1)$ to $Y_2 \setminus \text{nd}(K_2)$ along their
torus boundary, in such a way that the meridian $m_1$ of $K_1$ is
identified with the longitude $l_2$ of $K_2$ and vice versa. The
knot $K_1$ will induce a knot $K$ inside the three-manifold $Y$
obtained as above. We will denote the result of this construction
by $(Y_1,K_1)\perp (Y_2,K_2)=(Y,K)$. The following result
describes $\CFK(Y,K)$ and $\CFL(Y,K)$.

\begin{thm}
Suppose that $(Y,K)=(Y_1,K_1)\perp (Y_2,K_2)$ as above. Then
$$\RelSpinC(Y,K)\oplus \Z=\RelSpinC(Y_1,K_1) \oplus \RelSpinC(Y_2,K_2),$$
with $\Z$ being generated by $\text{PD}[m_2]$. Furthermore, for any
relative $\SpinC$ structure $\relspinc\in \RelSpinC(Y,K)$ we have
\begin{displaymath}
\begin{split}
\CFL(Y,K,\relspinc)&=\bigoplus_{\substack{\relspinc_i \in \RelSpinC(Y_i,K_i)\\
(\relspinc_1,\relspinc_2) \text{ induce } \relspinc}} \CFL(Y_1,K_1,\relspinc_1)
\otimes \CFK(Y_2,K_2,\relspinc_2),\\
\CFK(Y,K,\relspinc)&=\frac{\Big [\bigoplus_{\substack{\relspinc_i \in \RelSpinC(Y_i,K_i)\\
(\relspinc_1,\relspinc_2) \text{ induce } \relspinc}} \sC(\relspinc_1,\relspinc_2) \Big ]}
{\Large \sim},\\
\sC(\relspinc_1,\relspinc_2)&=\frac{\CFL(Y_1,K_1,\relspinc_1)\otimes
\CFK(Y_2,K_2,\relspinc_2)}{\CFL_0(Y_1,K_1,\relspinc_1)\otimes
\CFK_0(Y_2,K_2,\relspinc_2)},\\
\CFL_-(Y_1,K_1,\relspinc_1)&\otimes
\frac{\CFK(Y_2,K_2,\relspinc_2)}{\CFK_0(Y_2,K_2,\relspinc_2)}{\Large \sim}\\
&\frac{\CFL(Y_1,K_1,\relspinc_1)}{\CFL_0(Y_1,K_1,\relspinc_1)}\otimes
\CFK_-(Y_2,K_2,\relspinc_2-\text{PD}[m_2]).
\end{split}
\end{displaymath}
As before the equivalence relation $\sim$ is induced by the connecting isomorphisms
$\CFK_-\simeq \frac{\CFK}{\CFK_0}$ and $\CFL_-\simeq \frac{\CFL}{\CFL_0}$.
\end{thm}

\begin{proof}
If $H_1$ and $H_2$ are Heegaard diagrams for $(Y_1,K_1)$ and
$(Y_2,K_2)$ respectively, we may connect them by a one-handle
which is attached to the surfaces at the intersection of the
meridian and the longitude in each diagram. If we take the two
connected sums of curves $m_1\#l_2$ and $l_1\#m_2$ using this
one-handle, in such a way that they stay disjoint from each other,
the result will be a Heegaard diagram for $Y$. The meridian of $K$
in this Heegaard diagram may be thought of as the curve $m_1
\#l_2$, and the longitude may be represented by a parallel copy of
$l_1$ that stays on $\Sig_1$. It is now clear that for a
computation of $\CFL(Y,K)$ one can follow the argument of
Ozsv\'ath and Szab\'o, and that for a computation of $\CFK(Y,K)$,
the proof of theorem~\ref{thm:gluing} may be followed almost word
by word.
\end{proof}

\section{Examples; Alternating knots}

We have shown in the previous sections how our new $\{-,0,+\}$
filtrations of the complexes $\CFK(Y,K)$ and $\CFL(Y,K)$ play a
major role in gluing formulas. This role seems to be enough of a
justification for introducing them. However, the next question is
{\emph{how hard are these filtered complexes to compute?}}. The
computations of Ozsv\'ath and Szab\'o (\cite{OS-alternating}) and
Rasmussen (\cite{Ras,Ras2}) illustrated how well-behaved the
alternating knots, and more generally the perfect knots, are for
the purpose of Floer homology computations. So, it is natural to
turn to the computation of the $\{-,0,+\}$-filtered complexes
associated with
an alternating knot in $S^3$ (or more generally, perfect knots)
as the first computation.\\

In this section we will first compute the complexes $\CFL_\bullet(S^3,K)$ for an
alternating knot $(S^3,K)$:
\begin{thm}\label{thm:CFL(alternating)}
Suppose that $(S^3,K)$ is an alternating knot in $S^3$. Then the
filtered chain complex $\CFL_\bullet(S^3,K)$ is completely
determined by the Alexander polynomial $\Delta_K(t)$ and the
signature $\sigma(K)$.
\end{thm}

\begin{proof}
The proof is quite similar to our computation of
$\CFL(T(2,2n+1))$, where  $T(2,2n+1)$ is the $(2,2n+1)$ torus knot
(c.f. \cite{longitude}). Consider a standard Heegaard diagram
$$H_0=(\Sig',\alphas,\betas_0\cup\{\beta_1\},z)$$
for $(S^3,K)$, which is obtained from an alternating projection of
the knot $K$. The longitude $l$ of $(S^3,K)$ may be realized as a
curve which does not cut any of the curves in $\betas_0$, but it
cuts the meridian $m=\beta_1$ exactly once at the marked point
$z$. We may spin $l$ along the meridian $m$ for a large number $N$
of times to obtain a Heegaard diagram which is more appropriate
for our purposes. Now, add a new one-handle to the Heegaard
diagram so that an allowed Heegaard diagram  for the computation
of $\CFL(K)$ is obtained. There are two new curves added to the
Heegaard diagram, which will be denoted by $\lambda$ (the new
$\alpha$ curve) and $\mu$ (the new $\beta$ curve, obtained from
$m$). Thus we get the Heegaard diagram
$$H=(\Sig,\alphas \cup \{\lambda\}, \betas_0 \cup\{\mu\}\cup \{l\},z).$$
The final picture will
locally look like figure~\ref{fig:CFL(alternating)}.\\

\begin{figure}
\mbox{\vbox{\epsfbox{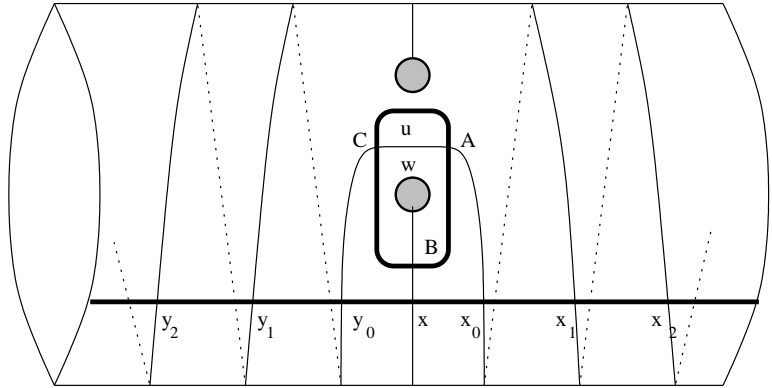}}}
\caption{\label{fig:CFL(alternating)} { The longitude of the
alternating knot $K$ may spin along the meridian several times to
generate the intersection points $x_0,x_1,...$ and $y_0,y_1,...$
with $\alpha_1$. Then by adding a handle, we may construct an
allowed Heegaard diagram for the computation of
$\CFL_\bullet(S^3,K)$. }}
\end{figure}

Let us denote the set of generators of the alternating projection
of $K$ by $\sZ$. The elements of $\sZ$ are in one-to-one
correspondence with the generators of the complex associated with
the initial Heegaard diagram $H_0$. In $H_0$ there is a unique
$\alpha$ circle which cuts $m$. We may assume that this curve is
named $\alpha_1$. Spinning $l$ along $m$ creates $2N$ new
intersections between $l$ and $\alpha_1$. Denote these
intersection points by $x_0,x_1,...,x_{N-1}$ on the right, and by
$y_0,y_1,...,y_{N-1}$ on the left. As usual, denote the
intersection points on $\lambda$ by $A,B$ and $C$, and
correspondingly divide the generators of the chain complex into
groups $\sA,\sB$ and $\sC$. Finally, denote the unique
intersection point of $\mu$ and $\alpha_1$ by $x$.\\

Spinning of $l$ along $m$ implies (c.f.
\cite{OS-3m2,OS-knot,Ras,longitude}) that the generators in the
relevant $\SpinC$ structures should contain one of the
intersection points $x,x_0,...,x_{N-1}$, or one of
$y_0,...,y_{N-1}$. Let us denote the relevant generators in
$\sA,\sB$ and $\sC$, by $\sA_0,\sB_0$ and $\sC_0$ respectively.
Then it is easy to check that we may make the identifications
\begin{displaymath}
\begin{split}
&\sA_0=\Big\{\{A,x,\z\} \ \Big |\ \  \z \in \sZ \Big\},\\
&\sC_0=\Big\{\{C,x,\z\} \ \Big |\ \  \z \in \sZ \Big\}, \text{ and}\\
&\sB_0=\Big\{\{B,x_i,\z\},\{B,y_i,\z\}\ \Big| \ \ \z \in \sZ,\ \
i=0,1,...,N-1\Big\}.
\end{split}
\end{displaymath}
Suppose that $F$ is a Seifert surface for the knot $K$, in $S^3$.
Such a surface may be obtained by capping off in $S^3\setminus
\text{nd}(K)$ a periodic domain $\sD$ for $K$ which has a zero
coefficient at $w$ and coefficient $1$ at $u$. We have $\sD=\sD_1
-\sD_2$, where $\sD_1$ is the small domain containing $u$.
 Then the
difference between the Maslov gradings of the two generators
$\{A,x,\z\}$ and $\{C,x,\z\}$ may be computed as follows. If we had illegally used
the domain $\sD_1$ to connect these two generators, the difference between the Maslov
gradings would have been one. In reality, we have to use the more complicated domain
$\sD_2$. Thus
$$\mu(\{A,x,\z\})-\mu(\{C,x,\z\})=
\mu(\sD_2)=1-\langle c_1(\spinc(\z)),F\rangle. $$ The number
$\langle c_1(\spinc(\z)),F\rangle$ is precisely  $2s(\z)$, where
$s(\z)$ is the $\Z$-grading induced by the knot $K$ on the
complex. Since $K$ is alternating and the diagram $H_0$ is
obtained from an alternating projection,
$s(\z)=\mu(\z)-\frac{\sigma(K)}{2}$, where $\sigma(K)$ denotes the
signature of the knot $K$ (see \cite{OS-alternating}). Here we are
already using the fact the the $\SpinC$ structure
$$\spinc(\{A,x,\z\})\in \RelSpinC(S^3,K)=\SpinC(S^3_0(K))$$ is the
same as the $\SpinC$ structure $\spinc(\z)$ assigned to $\z$ via
the
Heegaard diagram $H_0$.\\

In fact the $\SpinC$ structure assigned by the map $\spinc=\spinc_u$ to the above
generators may be computed via the formulas
\begin{displaymath}
\begin{split}
&\spinc( \{A,x,\z\})=\spinc(\z),\\
&\spinc( \{C,x,\z\})=\spinc(\z), \text{ and}\\
&\spinc( \{B,x_i,\z\})=
\spinc( \{B,y_i,\z\})= \spinc(\z)-i \text{PD}[m],\\
\end{split}
\end{displaymath}
where PD$[m]$ denotes the Poincar\'e dual of the meridian $m$ of the knot.\\

A consideration similar to that of section 3 of \cite{longitude}
shows that the Maslov grading levels of $\{B,x_i,\z\}$ and
$\{B,y_i,\z\}$ are independent of $i$. Also, because of the
trivial disks between $\{B,x_0,\z\}$ and $\{A,x,\z\}$, and between
$\{B,y_0,\z\}$ and $\{C,x,\z\}$ we know that their Maslov gradings
differ by one. We may summarize all these information as
\begin{displaymath}
\begin{split}
&\mu( \{A,x,\z\})=1+\sigma(K)-\mu(\z),\\
&\mu( \{C,x,\z\})=\mu(\z), \\
&\mu( \{B,x_i,\z\})=\sigma(K)-\mu(\z),\text{ and}\\
&\mu( \{B,y_i,\z\})=\mu(\z)-1.\\
\end{split}
\end{displaymath}
The boundary maps between the generators of $\sA_0$ and between the generators of
$\sC_0$ are exactly the boundary maps between the generators in $\sZ$, considered as
the generators of the complex associated with $H_0$. The boundary maps going from
$\sA_0$ to $\sB_0$ will give a map on the chain complex which is chain homotopic
to the map sending $\{A,x,\z\}$ to $\{B,x_0,\z\}$. Similarly, the map going from
$\sB_0$ to $\sC_0$ is chain homotopic to the map sending $\{B,y_0,\z\}$ to
$\{C,x,\z\}$. In order to complete the computation, we should study the maps within the
generators of type B, i.e. elements of $\sB_0$.\\

We remind the reader of the symmetry within the non-filtered
complex $\CFL(K)$, which was proved in \cite{longitude}. If we
re-consider the $\SpinC$ class of the generator $\x$ as the
average $\frac{1}{2}(\spinc_u(\x)+\spinc_w(\x))$ of the $\SpinC$
structures in $\RelSpinC(S^3,K)$ associated to $\x$ via the maps
$\spinc_u$ and $\spinc_w$, then we have shown in \cite{longitude}
that the non-filtered complexes $\CFL(S^3,K,\spinc)$ and
$\CFL(S^3,K,-\spinc)$ are in fact isomorphic, where $\spinc \in
\frac{1}{2}\PD[m]+\SpinC(S^3_0(K))$. Here $m$ is  the meridian of
$K$, and $-\spinc$ denotes the element with the property that
formally
$c_1(\spinc)+c_1(-\spinc)=0$.\\

Basically the same argument may be applied to show the symmetry of
the filtered version in the presence of the generators of
different types A,B, and C.  It is thus enough, as we did in the
computation of $\CFL(S^3,K)$ for the $(2,2n+1)$ torus knot
$K=T(2,2n+1)$ in $S^3$, to compute the complexes $\CFL(K,\spinc)$
for those elements $\spinc\in \frac{1}{2}\PD[m]+\SpinC(S^3_0(K))$
with $c_1(\spinc)>0$. Note that here we are taking the formal
first Chern class of an element of the form
$$\frac{1}{2}\PD[m]+\relspinc \in \frac{1}{2}\PD[m]+\RelSpinC(S^3,K)\simeq
\frac{1}{2}\PD[m]+\SpinC(S^3_0(K))$$ to be the cohomology class
with rational coefficients which evaluates on  a capped Seifert
surface $F$ of $K$ to give
$\frac{1}{2}+\langle c_1(\relspinc), F \rangle$.\\

In this regard, it is sufficient to compute the whole filtered complex for the
\emph{positive} $\SpinC$ classes $\spinc$ (i.e. those with
$\langle c_1(\spinc),F\rangle >0$).\\

When this is the case, there are no boundary maps going from the
generators of the form $\{B,x_i,\z\}$ to the generators of the
form $\{B,y_j,\w\}$. This is because the Maslov gradings of the
generators of the first type  are always less than the Maslov
gradings of the generators of the second type, as long as these
two generators are in the same $\SpinC$ class, and the associated
$\SpinC$ class $\spinc$ is a positive one. We refer
the reader to \cite{longitude} for a more detailed argument.\\

On the other hand, for orientational reasons, the domain of any
possible boundary map going from $\{B,y_j,\w\}$ to $\{B,x_i,\z\}$
has to have negative coefficients in the cylindrical region of
figure~\ref{fig:CFL(alternating)}, if both $i$ and $j$ are small,
compared to the number of twists (which may be chosen to be
arbitrarily large). The result of this discussion is that the only
boundary maps we have to understand for the positive $\SpinC$
structures $\spinc$, are the boundary maps within the generators
of the form $\{B,x_i,\z\}$ (which we call the generators of type
BX), as well as the boundary maps within the generators of the
form
$\{B,y_j,\w\}$ (which will be called generators of type BY).\\

\begin{figure}
\mbox{\vbox{\epsfbox{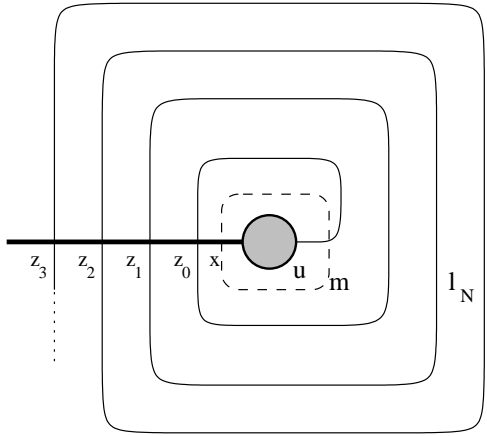}}}
\caption{\label{fig:alternating-spin} { If we redraw the diagram
$H_N$, around the marked point $u$, there will be the spined
longitude $l_N$ cutting $\alpha_1$ in several points
$z_0,z_1,...$, and cutting the meridian $m$ in $x$. $\alpha_1$ is
denoted by the bold curve, and $m$ by the dotted curve. }}
\end{figure}

Consider the Heegaard diagram obtained from $H_0$, by replacing
the meridian by the longitude $l_N=l+N.m$, which is obtained from
$l$ by spinning it $N$ times around the meridian $m$. Let us
denote this Heegaard diagram by $H_N$:
$$H_N=(\Sig',\alphas,\betas_0\cup\{l_N\},u).$$
Here $u$ denotes one of the marked points of the previous Heegaard
diagram $H_0$. Note that the Floer homology of this Heegaard
diagram computes $\CF(S^3_N(K))$,
of the three-manifold obtained from $S^3$ by a $N$-surgery on the knot $K$.\\
We may re-draw this Heegaard diagram around the marked point $u$, so that the
picture looks like figure~\ref{fig:alternating-spin}. The curve $l_N$ will be the
curve spinning around the center, and the unique $\alpha$-circle which cuts
$m$ in a single point, is the bold curve cutting $l_N$ in several (in fact $N$)
intersection points $z_0,z_1,...,z_{N-1}$. Again we assume that $z_0$ is the
closest to the center of the spiral.\\

If we consider a sufficiently large number $N$, then except for a fixed number
$k$ (independent of $N$) of the $\SpinC$ classes $\spinct \in \SpinC(S^3_N(K))$,
all the generators in the $\SpinC$ class $\spinct$ will contain an intersection
point $z_i$ on $l_N$. In fact, one of these $\SpinC$ classes will be generated
by the elements
$$\Big\{\{z_{g+s(\z)},\y\}\ \Big |\ \ \z=\{x,\y\}\in \sZ \Big\},$$
where $g$ is the genus of the alternating knot $K$. \\

Note that the generator $\{z_i,\y\}$ is quite similar to the
generator corresponding to $\{B,y_i,\z\}$ of type BY. In fact the
complex generated by the generators of type BY, and with the
boundary map induced from $\CFL(K)$, is isomorphic to the
complex generated by the generators of $\CF(S^3_N(K))$ of the above form.\\

We will first try to understand the part of $\CF(S^3_N(K))$ in the
$\SpinC$ class $\spinct$ considered above. Then we will use this
information to understand the boundary maps within the generators
of type BX and within the generators of type
BY.\\

Denote by $\betas_N$ the set of simple  closed curves
$\betas_0\cup\{l_N\}$, and let $\gammas$ be another set of simple
closed curves obtained as the union of  Hamiltonian isotopes
$\gammas_0$ of $\betas_0$, with the single curve $\{m\}$. Consider
the triple Heegaard diagram
$$\Hbar_N=(\Sig',\alphas,\gammas,\betas_N;u),$$
and the corresponding chain map
$$\sF:\CF(\alphas,\gammas)\otimes \CF(\gammas,\betas_N)\lra \CF(\alphas,\betas_N),$$
obtained by a count of holomorphic triangles associated with the Heegaard diagram
$\Hbar_N$. Note that the chain complex $\CF(\alphas,\gammas)$ is the filtered
chain complex $\CF(S^3)$, which is filtered using the alternating knot $K$.
Moreover the chain complex $\CF(\gammas,\betas_N)$ computes the Floer homology of the
three-manifold $\#^{g-1}(S^1\times S^2)$. Let $\Theta$ denote the top generator of this
homology group in the complex $\CF(\gammas,\betas_N)$. Moreover, denote the part of
$\sF(\x \otimes \Theta)$ in $\SpinC$ class $\spinct$ by $\sF_\spinct(\x)$. This way,
we get a chain map
$$\sF_\spinct:\CF(K,S^3) \lra \CF(S^3_N(K),\spinct),$$
where we use $\CF(K,S^3)$ to emphasize that our complex $\CF(S^3)$
admits a filtration
coming from the knot $K$.\\

It is quite straight forward to check, using the energy filtration
of \cite{OS-3m2}, that for any generator $\z=\{x,\y\}$ of the
complex $\CF(K,S^3)$ we have
\begin{equation}\label{eq:sFs-energy}
\sF_\spinct(\z)=\{z_{g+s(\z)},\y\}+\text{lower energy terms}.
\end{equation}
This is enough to show that the map $\sF_\spinct$ induces an
isomorphism of the chain complexes. Moreover, note that using the
isomorphism $\sF_\spinct$ and the filtration of $\CF(S^3)$ by $K$,
we may induce a $\Z$-filtration on the target complex
$\CF(S^3_N(K),\spinct)$.  Using  equation~\ref{eq:sFs-energy} it
is not hard to show that this filtration, and the filtration given
by
$$\text{filtration level of }\{z_i,\y\}=\text{ filtration level of }\{x,\y\},$$
give isomorphic filtered chain complexes.\\

Let us denote by $\sG^-(K,\ell)$ the subcomplex of $\CF(K,S^3)$
generated by the generators in filtration levels less than $\ell
\in \Z$. Denote the quotient complex by
$$\sG^+(K,\ell)=\frac{\CF(K,S^3)}{\sG^-(K,\ell)}.$$

For a positive $\SpinC$ structure $\spinc$, the generators of type BY which correspond
to $\spinc$ are precisely those $\{B,y_i,\z\}$ such that $\spinc(\z)-i\PD[m]=\spinc$. If
we use the isomorphisms
$$\Z\simeq \SpinC(S^3_0(K)),\text{ and}$$
$$\frac{1}{2}+\Z\simeq \frac{1}{2}\PD[m]+\SpinC(S^3_0(K)),$$
and denote the images of $\spinc(\z)$ and $\spinc$ in $\Z$ and $\frac{1}{2}+\Z$
under these
isomorphism by $s(\z)$ and $s$ respectively, then the relevant generators of type
BY are precisely
$$\Big\{\{B,y_{s(\z)-s-\frac{1}{2}},\z\}\ \Big|\ \ \z\in \Z,\ \ s(\z)> s\Big\}.$$
The corresponding complex generated using the generators of
$\CF(S^3_N(K),\spinct)$ will be $\sG^+(K,s+\frac{1}{2})$ and in
fact we see that this correspondence is an
isomorphism of the chain complexes.\\

Similarly, we may prove that the generators of type BX in the
$\SpinC$ class corresponding to $s\in \frac{1}{2}+\Z$ generate a
complex isomorphic to $\sG^-(K,-s+\frac{1}{2})$. Since we already
know the differentials between the generators of different types,
and all the complexes involved in this computation may be
determined using the filtered chain complex $\CFK(K,S^3)$, we may
argue as  in \cite{OS-alternating,Ras2} that the Alexander
polynomial $\Delta_K(t)$ of $K$, and its signature $\sigma(K)$,
determine the $\{-,0,+\}$-filtered chain complex
$\CFL_\bullet(S^3,K)$.
\end{proof}

In fact, if we continue to show the part of the $\Z$-filtered
complex $\CF(S^3)$ in filtration level less than $\ell$ by
$\sG^-(K,\ell)$ and its quotient complex by
$$\sG^+(K,\ell)=\frac{\CF(K,S^3)}{\sG^-(K,\ell)},$$
then we may explicitly determine the complex $\CFL_\bullet$ as is
stated in the following proposition.

\begin{prop}
Suppose that $s \in \frac{1}{2}+\Z$ represents a $\SpinC$ structure
$$\spinc \in \frac{1}{2}\PD[m]+\SpinC(S^3_0(K))\simeq \frac{1}{2}+\Z$$
for an alternating knot $K$, and suppose that the complexes $\sC_\bullet(s)$ for
$\bullet\in {-,0,+}$ are defined as
\begin{displaymath}
\begin{split}
\sC_-(s)=\frac{\sG^-(K,-s+\frac{1}{2})}{\sG^-(K,-s-\frac{1}{2})}&\\
\sC_0(s)=\Big [\frac{\sG^-(K,-s+\frac{1}{2})}{\sG^-(K,-s-\frac{1}{2})}\Big]&
\bigoplus \Big[\sG^-(K,-s+\frac{1}{2})\oplus \sG^+(K,s+\frac{1}{2})\Big]\\
\sC_+(s)=\Big [\frac{\sG^-(K,-s+\frac{1}{2})}{\sG^-(K,-s-\frac{1}{2})}\Big]&
\bigoplus \Big[\sG^-(K,-s+\frac{1}{2})\oplus \sG^+(K,s+\frac{1}{2})\Big]\\
&\bigoplus \Big [ \frac{\sG^-(K,s+\frac{3}{2})}{\sG^-(K,s+\frac{1}{2})}\Big].
\end{split}
\end{displaymath}
Then the complex $\CFL_\bullet(S^3,K,\spinc)$ is given by
$\sC_\bullet(s)$, with its boundary map modified by adding the
maps
\begin{displaymath}
\begin{split}
&d_{0-}:\sC_0 \lra \sC_-\\
&\sG^-(K,-s+\frac{1}{2})\xra{d_{0-}}\frac{\sG^-(K,-s+\frac{1}{2})}
{\sG^-(K,-s-\frac{1}{2})}, \text{ and}\\
&d_{+0}:\sC_+ \lra \sC_0\\
&\frac{\sG^-(K,s+\frac{3}{2})}{\sG^-(K,s+\frac{1}{2})} \xra{d_{+0}}
\sG^+(K,s+\frac{1}{2})=\frac{\CF(K,S^3)}{\sG^-(K,s+\frac{1}{2})}.
\end{split}
\end{displaymath}
Here $\bullet \in \{-,0,+\}$.
\end{prop}
\begin{proof}
The proof is just a simple algebraic argument, once we use our
information on the structure of the boundary maps between
different types of generators in positive $\SpinC$ structures,
which was obtained in  the proof of
theorem~\ref{thm:CFL(alternating)}. The symmetry with respect to
the $\SpinC$ structures implies the theorem in general.
\end{proof}

We may modify the Heegaard diagram $H_0$ in a different way to
obtain an allowed Heegaard diagram for the computation of
$\{-,0,+\}$-filtered complex $\CFK(S^3,K)$. Namely, we may spin
the longitude along the meridian several times, and then by adding
a one handle and a pair of curves to the diagram, change it to an
allowed Heegaard
diagram for the computation of $\CFK_\bullet(S^3,K)$.\\

An argument almost identical to the argument used for proving
theorem~\ref{thm:CFL(alternating)} shows the following:

\begin{thm}
Suppose that $(S^3,K)$ is an alternating knot in $S^3$, and let $\Delta_K(t)$ and
$\sigma(K)$ denote the Alexander polynomial and the signature of the knot $K$.
Then the homotopy type of the filtered chain complex $\CFK_\bullet(S^3,K)$ is
completely determined by $\Delta_K(t)$ and $\sigma(K)$.
\end{thm}

\section{Filtration of Floer homology and genus formulas}

This section is a short remark on how the genus of a knot $K$ is
related to the new $\{-,0,+\}$-filtration of the complexes
$\CFK(K)$ and $\CFL(K)$. In \cite{longitude} we used a technique
similar to the one used in the previous section, to obtain a sharp
genus bound from $\CFL(K)$. In fact, by spinning the longitude
along the meridian of a given knot $K$, we showed that the
non-triviality of the group $\CFK(K,g)$ implies the non-triviality
of $\CFK(K,g-\frac{1}{2})$, provided that $g>0$ is the genus of
$K$. The same argument may be used, almost word by word, to show
the following (slight) generalization, in presence of $\{-,0,+\}$
filtration.

For a complex $\sC$ which is graded by $\Z$ (in the sense that the
differential of $\sC$ takes elements in a grading level $s$ to
elements in grading level $s$), we define the degree $d_+(\sC)$ to
be given by
$$d_+(\sC)=\text{max}\Big\{s\in \Z \ |\ \ \sC(s)\text{ does not have trivial chain
homotopy type} \Big \},$$ and similarly define
$$d_-(\sC)=\text{min}\Big\{s\in \Z \ |\ \ \sC(s)\text{ does not have trivial chain
homotopy type} \Big \}.$$

\begin{thm}
Suppose that $K$ is a knot in $S^3$ of genus $g$. Then for the
complex $\CFK(K)=\CFK(S^3,K)$ we have
$$d_+(\CFK_-(K))+1=-d_-(\CFK_-(K))=g,$$
$$d_+(\CFK_0(K))=-d_-(\CFK_0(K))=g,\text{ and}$$
$$d_+(\CFK_+(K))=-d_-(\CFK_+(K))=g.$$
Furthermore, for the complex $\CFL(K)=\CFL(S^3,K)$ we have
$$d_+(\CFL_-(K))=-d_-(\CFL_-(K))=g,$$
$$d_+(\CFL_0(K))=-d_-(\CFL_0(K))=g,\text{ and}$$
$$d_+(\CFL_+(K))=1-d_-(\CFL_+(K))=g.$$
\end{thm}

We skip a more detailed proof of this last theorem, and refer the
reader to \cite{longitude}.


\begin{thebibliography}{99}


\bibitem{Burde} Burde, G., Zieschang, H., \emph{Knots}
de Gruyter Studies in Mathematics, 5.
Walter de Gruyter Co., Berlin, 2003.

\bibitem{longitude} Eftekhary, E., Longitude Floer homology and the Whitehead
double, \emph{preprint, available at} math.GT/0407211

\bibitem{FT} Fuchs, D., Tabachnikov, S.,
Invariants of Legandrian and transverse knots in the standard contact space,
\emph{Topology}
vol.36, No.5 (1997) 1025-1053

\bibitem{knot} Lickorish, W. B. R., \emph{An Introduction to knot theory},
Graduate texts in mathematics, 175, Springer-Verlag,New York, 1997

\bibitem{Liv} Livingston, C.,
Computations of the Ozsvath-Szabo knot concordance invariant,
\emph{preprint, available at} math.GT/0311036

\bibitem{OS-knot}  Ozsv\'{a}th, P., Szab\'o, Z.,
{Holomorphic disks and knot invariants},
\emph{to appear in Advances in Math., also available at} math.GT/0209056

\bibitem{OS-alternating}  Ozsv\'{a}th, P., Szab\'o, Z.,
Heegaard Floer homology and alternating knots,
\emph{Geometry and Topology}, Vol. 7 (2003), 225-254

\bibitem{OS-3m1}  Ozsv\'{a}th, P., Szab\'o, Z.,
{Holomorphic disks and topological invariants for closed three-manifolds},
\emph{to appear in Annals of Math., available at} math.SG/0101206

\bibitem{OS-3m2}  Ozsv\'{a}th, P., Szab\'o, Z.,
{Holomorphic disks and three-manifold invariants: properties and applications},
\emph{to appear in Annals of Math., available at} math.SG/0105202

\bibitem{OS-4m} Ozsv\'{a}th, P.,  Szab\'o, Z.,
{Holomorphic triangles and invariants for smooth four-manifolds}
\emph{ Duke Math. J.}  121  (2004),  no. 1, 1-34

\bibitem{OS-fibered} Ozsv\'{a}th, P., Szab\'o, Z.,
{Heegaard Floer homologies and contact structures},
\emph{preprint, available at} math.SG/0210127

\bibitem{OS-4genus} Ozsv\'{a}th, P., Szab\'o, Z.,
{Knot Floer homology and $4$-ball genus},
\emph{Geometry and Topology} vol.7, (2003) 615-639

\bibitem{OS-genus} Ozsv\'{a}th, P., Szab\'o, Z.,
{Holomorphic disks and genus bounds},\emph{ Geom. Topol.}
8  (2004), 311-334

\bibitem{Ras} Rasmussen, J.,
{Floer homology of surgeries on two bridge knots},
\emph{Alg. and Geo. Topology,2} (2002) 757-789

\bibitem{Ras2} Rasmussen, J.,
{Floer homology and knot complements},
\emph{Ph.D thesis, Harvard univ., also available at} math.GT/0306378

\bibitem{Rud} Rudolph, L.,
The slice genus and the Thurston-Bennequin invariant of a knot,
\emph{Proc. American Math. Soc.}, vol.125, No. 10 (1997) 3049-3050

\end{thebibliography}
\end{document}